\documentclass{amsart}%
\usepackage{amsfonts}
\usepackage{amsmath}
\usepackage{amssymb}
\usepackage{graphicx}%
\newtheorem{theorem}{Theorem}
\theoremstyle{plain}

\newtheorem{axiom}[theorem]{Assumption}

\newtheorem{corollary}[theorem]{Corollary}

\newtheorem{lemma}[theorem]{Lemma}

\newtheorem{proposition}[theorem]{Proposition}
\newtheorem{remark}[theorem]{Remark}

\newtheorem*{summary}{Proto-theorem}
\theoremstyle{definition}
\newtheorem*{acknowledgement}{Acknowledgments}

\newtheorem{example}{Example}
\numberwithin{equation}{section}

\begin{document}
\title[Local monotonicity and mean value formulas]{Local monotonicity and mean value formulas for evolving Riemannian manifolds}
\author{Klaus Ecker}
\address[Klaus Ecker]{ Freie Universit\"{a}t Berlin}
\email{ecker@math.fu-berlin.de}
\urladdr{http://geometricanalysis.mi.fu-berlin.de/people\_ecker.htm}
\author{Dan Knopf}
\address[Dan Knopf]{ University of Texas at Austin}
\email{danknopf@math.utexas.edu}
\urladdr{http://www.ma.utexas.edu/\symbol{126}danknopf/}
\author{Lei Ni}
\address[Lei Ni]{ University of California, San Diego}
\email{lni@math.ucsd.edu}
\urladdr{http://www.math.ucsd.edu/\symbol{126}lni/}
\author{Peter Topping}
\address[Peter Topping]{ University of Warwick}
\email{topping@maths.warwick.ac.uk}
\urladdr{http://www.maths.warwick.ac.uk/\symbol{126}topping/}

\begin{abstract}
We derive identities for general flows of Riemannian metrics that may be
regarded as local mean-value, monotonicity, or Lyapunov formulae. These
generalize previous work of the first author for mean curvature flow and other
nonlinear diffusions. Our results apply in particular to Ricci flow, where
they yield a local monotone quantity directly analogous to Perelman's reduced
volume $\tilde{V}$ and a local identity related to Perelman's average energy
$\mathcal{F}$.

\end{abstract}
\maketitle

\section{Introduction\label{Introduction}}

To motivate the local formulas we derive in this paper, consider the following
simple but quite general strategy for finding monotone quantities in geometric
flows, whose core idea is simply integration by parts. Let $(\mathcal{M}%
^{n},g(t))$ be a smooth one-parameter family of complete Riemannian manifolds
evolving for $t\in\lbrack a,b]$ by%
\begin{equation}
\frac{\partial}{\partial t}g=2h. \label{GeneralMetricEvolution}%
\end{equation}
Observe that the formal conjugate of the time-dependent heat operator
$\frac{\partial}{\partial t}-\Delta$ on the evolving manifold $(\mathcal{M}%
^{n},g(t))$ is $-(\frac{\partial}{\partial t}+\Delta+\operatorname*{tr}%
\!_{g}h)$. If $\varphi,\psi:\mathcal{M}^{n}\times\lbrack a,b]\rightarrow
\mathbb{R}$ are smooth functions for which the divergence theorem is valid
(e.g.~if $\mathcal{M}^{n}$ is compact or if $\varphi$ and $\psi$ and their
derivatives decay rapidly enough at infinity), one has\thinspace\footnote{Here
and throughout this paper, $d\mu$ denotes the volume form associated to
$g(t)$.}%
\begin{equation}
\frac{d}{dt}\int_{\mathcal{M}^{n}}\varphi\psi\,d\mu=\int_{\mathcal{M}^{n}%
}\{\psi\lbrack(\frac{\partial}{\partial t}-\Delta)\varphi]+\varphi
\lbrack(\frac{\partial}{\partial t}+\Delta+\operatorname*{tr}\!_{g}%
h)\psi]\}\,d\mu. \label{GlobalMonotone}%
\end{equation}
If $\varphi$ solves the heat equation and $\psi$ solves the adjoint heat
equation, it follows that the integral $\int_{\mathcal{M}^{n}}\varphi
\psi\,d\mu$ is independent of time. More generally, if $\psi\lbrack
(\frac{\partial}{\partial t}-\Delta)\varphi]$ and $\varphi\lbrack
(\frac{\partial}{\partial t}+\Delta+\operatorname*{tr}\!_{g}h)\psi]$ both have
the same sign, then $\int_{\mathcal{M}^{n}}\varphi\psi\,d\mu$ will be monotone
in $t$. If the product $\varphi\psi$ is geometrically meaningful, this can
yield useful results. Here are but a few examples.

\begin{example}
\label{Global-Heat}The simplest example uses the heat equation on Euclidean
space. Let
\begin{equation}
\psi(x,t)=\frac{1}{[4\pi(s-t)]^{n/2}}e^{-\frac{|y-x|^{2}}{4(s-t)}}\qquad
(x\in\mathbb{R}^{n},t<s) \label{BackwardHeatKernel}%
\end{equation}
denote the backward heat kernel with singularity at $(y,s)\in\mathbb{R}%
^{n}\times\mathbb{R}$. If $\varphi$ solves the heat equation and neither it
nor its derivatives grow too fast at infinity, then
\[
\varphi(y,s)=\lim_{t\nearrow s}\int_{\mathbb{R}^{n}}\varphi(x,t)\psi
(x,t)\,dx.
\]
Because $\frac{d}{dt}\int_{\mathbb{R}^{n}}\varphi(x,t)\psi(x,t)\,dx=0$, one
has $\varphi(y,s)=\int_{\mathbb{R}^{n}}\varphi(x,t)\psi(x,t)\,dx$ for all
$y\in\mathbb{R}^{n}$ and $t<s$, which illustrates the averaging property of
the heat operator.
\end{example}

\begin{example}
\label{Global-MCF}Let $F_{t}:\mathcal{M}^{n}\hookrightarrow\mathcal{M}_{t}%
^{n}\subset\mathbb{R}^{n+1}$ be a one-parameter family of hypersurfaces
evolving by mean curvature flow, $\frac{\partial}{\partial t}F_{t}=-H\nu$,
where $H$ is the mean curvature and $\nu$ the outward unit normal of the
hypersurface $\mathcal{M}_{t}^{n}$. This corresponds to $h=-HA$ in
(\ref{GeneralMetricEvolution}), where $A$ is the second fundamental form.
Define $\psi$ by formula~(\ref{BackwardHeatKernel}) applied to $x\in
\mathbb{R}^{n+1}$ and $t<s$. Using $\operatorname*{tr}\!_{g}h=-H^{2}$, one
calculates that%
\[
(\frac{\partial}{\partial t}+\Delta-H^{2})\psi=-\left\vert \frac{(x-y)^{\perp
}}{2\left(  s-t\right)  }-H\nu\right\vert ^{2}.
\]
Hence by (\ref{GlobalMonotone}),%
\[
\frac{d}{dt}\int_{\mathcal{M}_{t}^{n}}\varphi\psi\,d\mu=\int_{\mathcal{M}%
_{t}^{n}}[(\frac{\partial}{\partial t}-\Delta)\varphi]\psi\,d\mu
-\int_{\mathcal{M}_{t}^{n}}\left\vert \frac{(x-y)^{\perp}}{2\left(
s-t\right)  }-H\nu\right\vert ^{2}\varphi\psi\,d\mu.
\]
This is established for $\varphi\equiv1$ by Huisken \cite[Theorem
3.1]{Huisken90} and generalized by Huisken and the first author \cite[\S 1]%
{EH89} to any smooth $\varphi$ for which the integrals are finite and
integration by parts is permissible.

Hence $\int_{\mathcal{M}_{t}^{n}}\psi\,d\mu$ is monotone nonincreasing in time
and is constant precisely on homothetically shrinking solutions. The
monotonicity implies that the density $\Theta_{\mathcal{O}}^{\mathrm{MCF}%
}:=\lim_{t\nearrow0}\int_{\mathcal{M}_{t}^{n}}\psi\,d\mu$ of the limit point
$\mathcal{O}=(0,0)$ is well defined. Another consequence is that
$\sup_{\mathcal{M}_{b}^{n}}\varphi\leq\sup_{\mathcal{M}_{a}^{n}}\varphi$ if
$(\frac{\partial}{\partial t}-\Delta)\varphi\leq0$ for $t\in\lbrack a,b]$.
\end{example}

\begin{example}
\label{EntropyExample}A compact Riemannian manifold $(\mathcal{M}^{n},g(t))$
evolving by Ricci flow corresponds to $h=-\operatorname*{Rc}$ in
(\ref{GeneralMetricEvolution}), so that $\operatorname*{tr}\!_{g}h=-R$.
If\thinspace\footnote{Throughout this paper, $\nabla$ represents the
\emph{spatial }covariant derivative, and $\Delta=\operatorname*{tr}%
\!_{g}\nabla\nabla$.}
\[
\varphi\equiv1\qquad\text{and}\qquad\psi=\left[  \tau(2\Delta f-|\nabla
f|^{2}+R)+f-n\right]  (4\pi\tau)^{-n/2}e^{-f},
\]
then Perelman's entropy may be written as $\mathcal{W}\left(  g(t),f(t),\tau
(t)\right)  =\int_{\mathcal{M}^{n}}\varphi\psi\,d\mu$. If $d\tau/dt=-1$ and
$(\frac{\partial}{\partial t}+\Delta)f=|\nabla f|^{2}-R-\frac{n}{2\tau}$,
then
\[
(\frac{\partial}{\partial t}+\Delta-R)\psi=2|\operatorname*{Rc}+\nabla\nabla
f-\frac{1}{2\tau}g|^{2}(4\pi\tau)^{-n/2}e^{-f}.
\]
In this case, (\ref{GlobalMonotone}) becomes%
\[
\frac{d}{dt}\mathcal{W}(g(t),f(t),\tau(t))=\int_{\mathcal{M}^{n}%
}2|\operatorname*{Rc}+\nabla\nabla f-\frac{1}{2\tau}g|^{2}(4\pi\tau
)^{-n/2}e^{-f}\,d\mu,
\]
which is formula (3.4) of \cite{Perelman1}. In particular, $\mathcal{W}$ is
monotone increasing and is constant precisely on compact shrinking gradient solitons.
\end{example}

\begin{example}
\label{RVE}Again for $(\mathcal{M}^{n},g(t))$ evolving smoothly by Ricci flow
for $t\in\lbrack a,b]$, let $\ell$ denote Perelman's reduced distance
\cite{Perelman1} from an origin $(y,b)$. Take $\varphi\equiv1$ and choose
$\psi\equiv v$ to be the reduced-volume density\thinspace\footnote{The formula
used here and throughout this paper differs from Perelman's by the constant
factor $\left(  4\pi\right)  ^{-n/2}$. This normalization is more convenient
for our applications.}%
\[
v(x,t)=\frac{1}{[4\pi(b-t)]^{n/2}}e^{-\ell(x,b-t)}\qquad(x\in\mathcal{M}%
^{n},t<b).
\]
Then Perelman's reduced volume is given by $\tilde{V}\left(  t\right)
=\int_{\mathcal{M}^{n}}\varphi\psi\,d\mu$. By \cite[\S 7]{Perelman1},
$(\frac{\partial}{\partial t}+\Delta-R)v\geq0$ holds in the barrier sense,
hence in the distributional sense.\footnote{It is a standard fact that a
suitable barrier inequality implies a distributional inequality. See
\cite{RFV2} for relevant definitions and a proof. A direct proof for $v$ is
found in \cite[Lemma 1.12]{Ye}.} Thus one obtains monotonicity of the reduced
volume if $\mathcal{M}^{n}$ is compact or if its Ricci curvature is bounded.

More generally, one gets monotonicity of $\int_{\mathcal{M}^{n}}\varphi
\psi\,d\mu$ for any nonnegative supersolution $\varphi$ of the heat equation.
In particular, taking $\varphi(x,t)=R(x,t)-R_{\min}(0)$ on a compact manifold
and noting that $(\frac{\partial}{\partial t}-\Delta)\varphi\geq0$ holds
pointwise, one verifies that $\int_{\mathcal{M}^{n}}[R-R_{\min}(0)]v\,d\mu$ is
nondecreasing in time.

In \cite{FIN05}, Feldman, Ilmanen, and the third author introduce an expanding
entropy and a forward reduced volume for compact manifolds evolving by Ricci
flow. Monotonicity of these quantities may also be derived from
(\ref{GlobalMonotone}) with $\varphi\equiv1$.

Similar ideas play important roles in Perelman's proofs of differential
Harnack estimates \cite[\S 9]{Perelman1} and pseudolocality \cite[\S 10]%
{Perelman1}.
\end{example}

\bigskip

The strategy of integration by parts can be adapted to yield \emph{local
}monotone quantities for geometric flows. We shall present a rigorous
derivation in Section~\ref{Derivation} when we prove our main result,
Theorem~\ref{MainTheorem}. Before doing so, however, we will explain the
underlying motivations by a purely formal argument. Suppose for the purposes
of this argument that $\Omega=%
{\textstyle\bigcup\nolimits_{a\leq t\leq b}}
\Omega_{t}$ is a smooth, precompact subset of $\mathcal{M}^{n}\times\lbrack
a,b]$. Assume that $\partial\Omega_{t}$ is smooth with outward unit normal
$\nu$, and let $d\sigma$ denote the measure on $\partial\Omega_{t}$ induced by
$g(t)$. If the product $\varphi\psi$ vanishes on $\partial\Omega$, then%
\begin{align}
\int_{\Omega}  &  \{\psi\lbrack(\frac{\partial}{\partial t}-\Delta
)\varphi]+\varphi\lbrack(\frac{\partial}{\partial t}+\Delta+\operatorname*{tr}%
\!_{g}h)\psi]\}\,d\mu\,dt\label{STG2-1}\\
&  =\int_{a}^{b}\left(  \frac{d}{dt}\int_{\Omega_{t}}\varphi\psi\,d\mu\right)
\,dt+\int_{\partial\Omega}(\varphi\left\langle \nabla\psi,\nu\right\rangle
-\psi\left\langle \nabla\varphi,\nu\right\rangle )\,d\sigma\,dt.\nonumber
\end{align}
This formula may be regarded as a space-time analog of Green's second
identity. In the special case that $\Omega$ is the super-level set
$\{(x,t):\psi(x,t)>0\}$ and both $\Omega_{a}$ and $\Omega_{b}$ are empty, then
$\nu=-|\nabla\psi|^{-1}\nabla\psi$, whence (\ref{STG2-1}) reduces to%
\begin{equation}
\int_{\{\psi>0\}}\{\psi\lbrack(\frac{\partial}{\partial t}-\Delta
)\varphi]+\varphi\lbrack(\frac{\partial}{\partial t}+\Delta+\operatorname*{tr}%
\!_{g}h)\psi]\}\,d\mu\,dt+\int_{\{\psi=0\}}\varphi|\nabla\psi|\,d\sigma\,dt=0.
\label{STG2-2}%
\end{equation}
Formula (\ref{STG2-2}) enables a strategy for the construction of local
monotone quantities.

Here is the strategy, again presented as a purely formal argument. Let
$\varphi$ and $\Psi>0$ be given. Define $\psi=\log\Psi$, and for $r>0$, let
$\psi_{(r)}=\log(r^{n}\Psi)$. Take $\Omega$ to be the set $E_{r}$ defined for
$r>0$ by%
\begin{equation}
E_{r}:=\{(x,t):\Psi(x,t)>r^{-n}\}=\{(x,t):\psi_{(r)}>0\}. \label{Er-Formal}%
\end{equation}
(When $\Psi$ is a fundamental solution\thinspace\footnote{See
Section~\ref{MVT} below.} of a backward heat equation, the set $E_{r}$ is
often called a `heatball'.) Assume for the sake of this formal argument that
the outward unit normal to the time slice $E_{r}(t):=E_{r}\cap(\mathcal{M}%
^{n}\times\{t\})$ is $\nu=-|\nabla\psi|^{-1}\nabla\psi$. Observing that%
\begin{equation}
(\frac{\partial}{\partial t}+\Delta)\psi=\Psi^{-1}(\frac{\partial}{\partial
t}+\Delta+\operatorname*{tr}\!_{g}h)\Psi-|\nabla\psi|^{2}-\operatorname*{tr}%
\!_{g}h \label{log-eqn}%
\end{equation}
and applying the coarea formula to each time slice $E_{r}(t)$, followed by an
integration in $t$, one obtains%
\begin{equation}
\frac{d}{dr}\int_{E_{r}}|\nabla\psi|^{2}\varphi\,d\mu\,dt=\frac{n}{r}%
\int_{\partial E_{r}}|\nabla\psi|\varphi\,d\sigma\,dt. \label{coarea1}%
\end{equation}
Similarly, one has
\begin{subequations}
\label{coarea2}%
\begin{align}
\frac{d}{dr}\int_{E_{r}}(\operatorname*{tr}\!_{g}h)\psi_{(r)}\varphi
\,d\mu\,dt  &  =\int_{E_{r}}(\operatorname*{tr}\!_{g}h)[\frac{\partial
}{\partial r}\log(r^{n}\Psi)]\varphi\,d\mu\,dt\\
&  =\frac{n}{r}\int_{E_{r}}(\operatorname*{tr}\!_{g}h)\varphi\,d\mu\,dt,
\end{align}
because the boundary integral vanishes in this case. Now by rearranging
(\ref{STG2-2}) and using (\ref{log-eqn})--(\ref{coarea2}), one gets%
\end{subequations}
\begin{align*}
\int_{E_{r}}  &  \{\psi_{(r)}[(\frac{\partial}{\partial t}-\Delta
)\varphi]+\Psi^{-1}[(\frac{\partial}{\partial t}+\Delta+\operatorname*{tr}%
\!_{g}h)\Psi]\varphi\}\,d\mu\,dt\\
&  =\int_{E_{r}}\left\{  |\nabla\psi|^{2}-(\operatorname*{tr}\!_{g}%
h)\psi_{(r)}\right\}  \varphi\,d\mu\,dt-\int_{\partial E_{r}}\varphi
|\nabla\psi|\,d\sigma\,dt+\int_{E_{r}}(\operatorname*{tr}\!_{g}h)\varphi
\,d\mu\,dt\\
&  =\int_{E_{r}}\left\{  |\nabla\psi|^{2}-(\operatorname*{tr}\!_{g}%
h)\psi_{(r)}\right\}  \varphi\,d\mu\,dt-\frac{r}{n}\frac{d}{dr}\int_{E_{r}%
}\left\{  |\nabla\psi|^{2}-(\operatorname*{tr}\!_{g}h)\psi_{(r)}\right\}
\varphi\,d\mu\,dt.
\end{align*}
Defining%
\begin{equation}
P_{\varphi,\Psi}(r):=\int_{E_{r}}\left\{  |\nabla\log\Psi|^{2}%
-(\operatorname*{tr}\!_{g}h)\log(r^{n}\Psi)\right\}  \varphi\,d\mu\,dt
\label{DefineP-Formal}%
\end{equation}
and applying an integrating factor, one obtains the following formal identity.
Since $\log(r^{n}\Psi)=\psi_{(r)}>0$ in $E_{r}$, this identity produces a
local monotone quantity whenever $(\frac{\partial}{\partial t}-\Delta)\varphi$
and $\varphi(\frac{\partial}{\partial t}+\Delta+\operatorname*{tr}\!_{g}%
h)\Psi$ have the same sign.

\begin{summary}
Whenever the steps above can be rigorously justified and all integrals in
sight make sense, the identity
\begin{equation}
\frac{d}{dr}\left(  \frac{P_{\varphi,\Psi}(r)}{r^{n}}\right)  =-\frac
{n}{r^{n+1}}\int_{E_{r}}\left\{  \log(r^{n}\Psi)(\frac{\partial}{\partial
t}-\Delta)\varphi+\Psi^{-1}[(\frac{\partial}{\partial t}+\Delta
+\operatorname*{tr}\!_{g}h)\Psi]\varphi\right\}  \,d\mu\,dt
\label{ProtoTheorem}%
\end{equation}
will hold in an appropriate sense.
\end{summary}

In spirit, (\ref{ProtoTheorem}) is a parabolic analogue of the formula%
\[
\frac{d}{dr}\left(  \frac{1}{r^{n}}\int_{|x-y|<r}\varphi(x)\,d\mu\right)
=\frac{1}{2r^{n+1}}\int_{|x-y|<r}(r^{2}-|x-y|^{2})\Delta\varphi\,d\mu,
\]
which for harmonic $\varphi$ (i.e.~$\Delta\varphi=0$) leads to the classical
local mean-value representation formulae%
\[
\varphi(y)=\frac{1}{\omega_{n}r^{n}}\int_{|x-y|<r}\varphi(x)\,d\mu=\frac
{1}{n\omega_{n}r^{n-1}}\int_{|x-y|=r}\varphi(x)\,d\sigma.
\]

\bigskip

The main result of this paper, Theorem~\ref{MainTheorem}, is a rigorous
version of the motivational proto-theorem above. We establish
Theorem~\ref{MainTheorem} in a sufficiently robust framework to provide new
proofs of some classical mean-value formulae (Examples~\ref{Euclidean}%
--\ref{NiEntropies}), to generate several new results
(Corollaries~\ref{ReducedDistanceCorollary}, \ref{VolumeEquality},
\ref{SolitonReducedVolume}, \ref{FRD}, \ref{Entropy}, \ref{Moving},
\ref{LN-MVI}) and to permit generalizations for future applications. Our
immediate original results are organized as follows: in Section ~\ref{RV}, we
study Perelman's reduced volume for manifolds evolving by Ricci flow; in
Section~\ref{AE}, we consider Perelman's average energy for manifolds evolving
by Ricci flow; and in Section~\ref{MVT}, we discuss heat kernels on evolving
Riemannian manifolds (including fixed manifolds as an interesting special
case). In \cite{Ni-MVT}, the third author applies some of these results to
obtain local regularity theorems for Ricci flow. Potential future
generalizations that we have in mind concern varifold (Brakke) solutions of
mean curvature flow, solutions of Ricci flow with surgery, and fundamental
solutions in the context of `weak' (Bakry--\'{E}mery) Ricci curvature,
e.g.~\cite{LV06}.

\bigskip

As noted above, Theorem~\ref{MainTheorem} allows new proofs of several
previously known local monotonicity formulae, all of which should be compared
with (\ref{ProtoTheorem}). To wit:

\begin{example}
\label{Euclidean}Consider the Euclidean metric on $\mathcal{M}^{n}%
=\mathbb{R}^{n}$ with $h=0$. If $\Psi$ is the backwards heat kernel
(\ref{BackwardHeatKernel}) centered at $(y,s)$ and the heatball $E_{r}\equiv
E_{r}(y,s)$ is defined by (\ref{Er-Formal}), then (\ref{DefineP-Formal})
becomes%
\[
P_{\varphi,\Psi}(r)=\int_{E_{r}(y,s)}\varphi(x,t)\frac{|y-x|^{2}}{4(s-t)^{2}%
}\,d\mu\,dt.
\]
Thus (\ref{ProtoTheorem}) reduces to%
\begin{equation}
\frac{d}{dr}\left(  \frac{P_{\varphi,\Psi}(r)}{r^{n}}\right)  =-\frac
{n}{r^{n+1}}\int_{E_{r}(y,s)}\log(r^{n}\Psi)(\frac{\partial}{\partial
t}-\Delta)\varphi\,d\mu\,dt. \label{ProtoHeat}%
\end{equation}
Since $\int_{E_{r}(y,s)}\frac{|y-s|^{2}}{4(s-t)^{2}}\,d\mu\,dt=1$, this
implies the mean value identity%
\begin{equation}
\varphi(y,s)=\frac{1}{r^{n}}\int_{E_{r}(y,s)}\varphi(x,t)\frac{|y-x|^{2}%
}{4(s-t)^{2}}\,d\mu\,dt \label{MVI-heat}%
\end{equation}
for all $\varphi$ satisfying $(\frac{\partial}{\partial t}-\Delta)\varphi=0$.
This localizes Example~\ref{Global-Heat}.

To our knowledge, Pini \cite{Pini51, Pini54a, Pini54b} was the first to prove
(\ref{MVI-heat}) in the case $n=1$. This was later generalized to $n>1$ by
Watson \cite{Watson73}. The general formula (\ref{ProtoHeat}) appears in
Evans--Gariepy \cite{EG82}. There are many similar mean-value representation
formulae for more general parabolic operators. For example, see
Fabes--Garofalo \cite{FG87} and Garofalo--Lanconelli \cite{GL88}. (Also see
Corollaries~\ref{Moving} and \ref{LN-MVI}, below.)
\end{example}

\begin{example}
Surface integrals over heatballs first appear in the work of Fulks
\cite{Fulks66}, who proves that a continuous function $\varphi$ on
$\mathbb{R}^{n}\times(a,b)$ satisfies%
\[
\varphi(y,s)=\frac{1}{r^{n}}\int_{\partial E_{r}}\varphi(x,t)\frac{|y-x|^{2}%
}{\sqrt{4|y-x|^{2}(s-t)^{2}+[|y-x|^{2}-2n(s-t)]^{2}}}\,d\sigma
\]
for all sufficiently small $r>0$ if and only if $\varphi$ is a solution of the
heat equation. (Compare to Corollary~\ref{LeiNi2} below.)
\end{example}

\begin{example}
\label{MCFdensity}Previous results of the first author \cite{Ecker01} localize
Example~\ref{Global-MCF} for mean curvature flow. On $\mathbb{R}^{n+1}%
\times(-\infty,0)$, define $\Psi(x,t):=(-4\pi t)^{-n/2}e^{|x|^{2}/4t}$.
Substitute $\Psi^{-1}(\frac{\partial}{\partial t}+\Delta+\operatorname*{tr}%
\!_{g}h)\Psi=-|\nabla^{\perp}\psi+H\nu|^{2}$ and $\operatorname*{tr}%
\!_{g}h=-H^{2}$ into (\ref{DefineP-Formal}) and (\ref{ProtoTheorem}). If the
space-time track $\mathcal{M}=%
{\textstyle\bigcup\nolimits_{t<0}}
\mathcal{M}_{t}^{n}$ of a solution to mean curvature flow is well defined in
the cylinder $B(0,\sqrt{2n\bar{r}^{2}/\pi})\times(-\bar{r}^{2}/4\pi,0)$, then
\cite{Ecker01} proves that formula~(\ref{ProtoTheorem}), with the integrals
taken over $E_{r}\cap\mathcal{M}$, holds in the distributional sense for any
$r\in(0,\bar{r})$ and any $\varphi$ for which all integral expressions are
finite. In particular, $P_{1,\Psi}(r)/r^{n}$ is monotone increasing in $r$.
The density $\Theta_{\mathcal{O}}^{\mathrm{MCF}}:=\lim_{t\nearrow0}%
\int_{\mathcal{M}_{t}^{n}}\Psi(x,t)\,d\mu$ of the limit point $\mathcal{O}%
=(0,0)$ can thus be calculated locally by
\[
\Theta_{\mathcal{O}}^{\mathrm{MCF}}=\lim_{r\searrow0}\frac{P_{1,\Psi}%
(r)}{r^{n}}.
\]
(Compare to Corollary~\ref{SolitonReducedVolume} below.)

Related work of the first author for other nonlinear diffusions is found in
\cite{Ecker05}.
\end{example}

\begin{example}
\label{NiEntropies}Perelman's scaled entropy $\mathcal{W}$ and the forward
reduced volume $\theta_{+}$ are localized by the third author
\cite[Propositions 5.2, 5.3, 5.4]{Ni-LYH}. Although only stated there for
K\"{a}hler--Ricci flow, these localizations remain valid for Ricci flow in
general. They are motivated by the first author's work on mean curvature flow
\cite{E2} and arise from (\ref{GlobalMonotone}) by taking $\varphi$ to be a
suitable cutoff function defined with respect to $\tau\ell$ and $\tau\ell_{+}%
$, respectively.
\end{example}

\bigskip

The remainder of this paper is organized as follows. In
Section~\ref{Derivation}, we rigorously derive Theorem~\ref{MainTheorem}: the
general local monotonicity formula motivated by formula~(\ref{ProtoTheorem})
above. In Section~\ref{LGE}, we derive a local gradient estimate for solutions
of the conjugate heat equation. In Sections~\ref{RV}--\ref{MVT}, we apply this
machinery to obtain new results in some special cases where our assumptions
can be checked and in which (\ref{ProtoTheorem}) simplifies and becomes more
familiar. The Appendix (Section \ref{Appendix}) reviews some relevant
properties of Perelman's reduced distance and volume.

\begin{acknowledgement}
K.E.~ was partially supported by SFB 647. D.K.~was partially supported by NSF
grants DMS-0511184, DMS-0505920, and a University of Texas Summer Research
Assignment. L.N.~ was partially supported by NSF grants and an Alfred P.~Sloan
Fellowship. P.T.~was partially supported by an EPSRC Advanced Research Fellowship.

L.N.~ thanks both Professor Bennett Chow and Professor Peter Ebenfelt for
bringing Watson's mean-value equality to his attention. This motivated him to
study heatball constructions and in particular \cite{FG87} and \cite{Ecker01}.
He also thanks Professor Peter Li for many helpful discussions.
\end{acknowledgement}

\section{The rigorous derivation\label{Derivation}}

Let $-\infty<a<b<\infty$, and let $(\mathcal{M}^{n},g(t))$ be a smooth
one-parameter family of complete Riemannian manifolds evolving by
(\ref{GeneralMetricEvolution}) for $t\in\lbrack a,b]$. As noted above, the
formal conjugate of the heat operator $\frac{\partial}{\partial t}-\Delta$ on
$(\mathcal{M}^{n},g(t))$ is $-(\frac{\partial}{\partial t}+\Delta
+\operatorname*{tr}\!_{g}h)$. For $\alpha\in\mathbb{R}$, we adopt the standard
notation $[\alpha]_{+}:=\max\{\alpha,0\}$.

Let $\Psi$ be a given positive function on $\mathcal{M}^{n}\times\lbrack
a,b)$. As in Section~\ref{Introduction}, it is convenient to work with
\begin{equation}
\psi:=\log\Psi
\end{equation}
and the function defined for each $r>0$ by%
\begin{equation}
\psi_{(r)}:=\psi+n\log r.
\end{equation}
For $r>0$, we define the space-time super-level set (`heatball')
\begin{subequations}
\label{DefineHeatBall}%
\begin{align}
E_{r}  &  =\{(x,t)\in\mathcal{M}^{n}\times\lbrack a,b):\Psi>r^{-n}\}\\
&  =\{(x,t)\in\mathcal{M}^{n}\times\lbrack a,b):\psi_{(r)}>0\}.
\end{align}

We would like to allow $\Psi$ to blow up as we approach time $t=b$; in
particular, we have in mind various functions which have a singularity that
agrees asymptotically with a (backwards) heat kernel centered at some point in
$\mathcal{M}^{n}$ at time $t=b$. (See Sections~\ref{AE}--\ref{MVT}.) In this
context, we make, for the moment, the following three assumptions about $\Psi$.
\end{subequations}
\begin{axiom}
\label{LL}The function $\Psi$ is locally Lipschitz on $\mathcal{M}^{n}%
\times\lbrack a,s]$ for any $s\in(a,b)$.
\end{axiom}

\begin{axiom}
\label{Omega}There exists a compact subset $\Omega\subseteq\mathcal{M}^{n}$
such that $\Psi$ is bounded outside $\Omega\times\lbrack a,b)$.
\end{axiom}

\begin{axiom}
\label{r-bar}There exists $\bar{r}>0$ such that%
\[
\lim_{s\nearrow b}\left(  \int_{E_{\bar{r}}\cap(\mathcal{M}^{n}\times
\{s\})}|\psi|\,d\mu\right)  =0
\]
and%
\[
\int_{E_{\bar{r}}}|\nabla\psi|^{2}\,d\mu\,dt<\infty.
\]

\end{axiom}

\begin{remark}
\label{CompactHeatballs}By the continuity of $\Psi$ from Assumption~\ref{LL}
and its boundedness from Assumption~\ref{Omega}, we can be sure, after
reducing $\bar{r}>0$ if necessary, that $\Psi\leq\bar{r}^{-n}$ outside some
compact subset of $\mathcal{M}^{n}\times(a,b]$. In particular, we then know
that the super-level sets $E_{r}$ lie inside this compact subset for
$r\in(0,\bar{r}]$.
\end{remark}

\begin{remark}
By Assumption~\ref{r-bar} and compactness of $E_{\bar{r}}$, one has
$\int_{E_{\bar{r}}}|\psi|\,d\mu\,dt<\infty$.
\end{remark}

\begin{remark}
We make no direct assumptions about the regularity of the sets $E_{r}$ themselves.
\end{remark}

Let $\varphi$ be an arbitrary smooth function on $\mathcal{M}^{n}\times(a,b]$.
By Assumption~\ref{r-bar}, the quantity
\begin{equation}
P_{\varphi,\Psi}(r):=\int_{E_{r}}[|\nabla\psi|^{2}-\psi_{(r)}%
(\operatorname*{tr}\!_{g}h)]\varphi\,d\mu\,dt \label{DefineP}%
\end{equation}
is finite for $r\in(0,\bar{r}]$. Our main result is as follows:

\begin{theorem}
\label{MainTheorem}Suppose that $(\mathcal{M}^{n},g(t))$ is a smooth
one-parameter family of complete Riemannian manifolds evolving by
(\ref{GeneralMetricEvolution}) for $t\in\lbrack a,b]$, that $\Psi
:\mathcal{M}^{n}\times\lbrack a,b)\rightarrow(0,\infty)$ satisfies
Assumptions~1--3, that $\bar{r}>0$ is chosen according to
Assumption~\ref{r-bar} and Remark~\ref{CompactHeatballs}, and that
$0<r_{0}<r_{1}\leq\bar{r}$.

If $\Psi$ is smooth and the function
\[
\frac{(\frac{\partial}{\partial t}+\Delta+\operatorname*{tr}\!_{g}h)\Psi}%
{\Psi}\equiv\frac{\partial\psi}{\partial t}+\Delta\psi+|\nabla\psi
|^{2}+\operatorname*{tr}\!_{g}h
\]
belongs to $L^{1}(E_{\bar{r}})$, then%
\begin{align}
\frac{P_{\varphi,\Psi}(r_{1})}{r_{1}^{n}}-\frac{P_{\varphi,\Psi}(r_{0})}%
{r_{0}^{n}}  &  =\label{Pexp1}\\
\int_{r_{0}}^{r_{1}}\frac{n}{r^{n+1}}\int_{E_{r}}[-(\frac{\partial\psi
}{\partial t}  &  +\Delta\psi+|\nabla\psi|^{2}+\operatorname*{tr}%
\!_{g}h)\varphi-(\psi+n\log r)(\frac{\partial\varphi}{\partial t}%
-\Delta\varphi)]\,d\mu\,dt\,dr.\nonumber
\end{align}

If, instead, $\Psi$ is merely locally Lipschitz in the sense of
Assumption~\ref{LL}, and the inequality
\[
\frac{(\frac{\partial}{\partial t}+\Delta+\operatorname*{tr}\!_{g}h)\Psi}%
{\Psi}\geq0
\]
holds in the distributional sense, and $\varphi\geq0$, then%
\begin{equation}
\frac{P_{\varphi,\Psi}(r_{1})}{r_{1}^{n}}-\frac{P_{\varphi,\Psi}(r_{0})}%
{r_{0}^{n}}\leq-\int_{r_{0}}^{r_{1}}\frac{n}{r^{n+1}}\int_{E_{r}}(\psi+n\log
r)(\frac{\partial\varphi}{\partial t}-\Delta\varphi)\,d\mu\,dt\,dr.
\label{Pexp2}%
\end{equation}

\end{theorem}

\begin{remark}
If $\varphi$ solves the heat equation and $\Psi$ solves the conjugate heat
equation, then (\ref{Pexp1}) implies that $P_{\varphi,\Psi}(r)/r^{n}$ is
independent of $r$. See Example~\ref{Euclidean} (above) and
Corollary~\ref{Moving} (below).
\end{remark}

\begin{proof}
[Proof of Theorem \ref{MainTheorem}]We begin by assuming that $\Psi$ is
smooth. In the proof, we write $P(\cdot)\equiv P_{\varphi,\Psi}(\cdot)$. For
most of the proof, we will work with a modified function, namely%
\begin{equation}
P(r,s):=\int_{E_{r}\cap(\mathcal{M}^{n}\times\lbrack a,s])}[|\nabla\psi
|^{2}-\psi_{(r)}(\operatorname*{tr}\!_{g}h)]\varphi\,d\mu\,dt,
\end{equation}
arising from restriction to the time interval $[a,s]$, for some $s\in(a,b)$.
As a result, we will only be working on domains on which $\Psi$ and its
derivatives are bounded, and the convergence of integrals will not be in
doubt. A limit $s\nearrow b$ will be taken at the end.

Let $\zeta:\mathbb{R}\rightarrow\lbrack0,1]$ be a smooth function with the
properties that $\zeta(y)=0$ for $y\leq0$ and $\zeta^{\prime}(y)\geq0$. Let
$Z:\mathbb{R}\rightarrow\lbrack0,\infty)$ denote the primitive of $\zeta$
defined by $Z(y)=\int_{-\infty}^{y}\zeta(x)\,dx$. One should keep in mind that
$\zeta$ can be made very close to the Heaviside function, in which case $Z(y)$
will lie a little below $[y]_{+}$.

For $r\in(0,\bar{r}]$ and $s\in(a,b)$, we define%
\begin{equation}
Q(r,s):=\int_{\mathcal{M}^{n}\times\lbrack a,s]}[|\nabla\psi|^{2}\zeta
(\psi_{(r)})-Z(\psi_{(r)})(\operatorname*{tr}\!_{g}h)]\varphi\,d\mu\,dt,
\label{DefineQ}%
\end{equation}
which should be regarded as a perturbation of $P(r,s)$, and will relieve us of
some technical problems arising from the fact that we have no control on the
regularity of $E_{r}$. Note that $\zeta(\psi_{(r)})$ and $Z(\psi_{(r)})$ have
support in $E_{r}$. Therefore, the convergence of the integrals is guaranteed.

In the following computations, we suppress the dependence of $Q$ on $s$ and
assume that each integral is over the space-time region $\mathcal{M}^{n}%
\times\lbrack a,s]$ unless otherwise stated. One has%
\begin{align}
\frac{r^{n+1}}{n}\frac{d}{dr}\left[  \frac{Q(r)}{r^{n}}\right]   &  =\frac
{r}{n}Q^{\prime}(r)-Q(r)\label{Qderiv1}\\
&  =\int[|\nabla\psi|^{2}\zeta^{\prime}(\psi_{(r)})-Z^{\prime}(\psi
_{(r)})(\operatorname*{tr}\!_{g}h)]\varphi\,d\mu\,dt-Q(r)\nonumber\\
&  =\int[|\nabla\psi|^{2}\zeta^{\prime}(\psi_{(r)})]\varphi\,d\mu
\,dt-\int[\zeta(\psi_{(r)})(\operatorname*{tr}\!_{g}h)]\varphi\,d\mu
\,dt\nonumber\\
&  -\int[|\nabla\psi|^{2}\zeta(\psi_{(r)})]\varphi\,d\mu\,dt+\int[Z(\psi
_{(r)})(\operatorname*{tr}\!_{g}h)]\varphi\,d\mu\,dt.\nonumber
\end{align}
The first integral and the last integral in the last equality on the
right-hand side require further attention.

For the first of these, we keep in mind that $\nabla\psi=\nabla\psi_{(r)}$ and
compute%
\begin{align}
\int &  [|\nabla\psi|^{2}\zeta^{\prime}(\psi_{(r)})]\varphi\,d\mu
\,dt\label{term1calc}\\
&  =\int\left\langle \nabla\psi,\nabla(\zeta(\psi_{(r)}))\right\rangle
\varphi\,d\mu\,dt\nonumber\\
&  =-\int[(\Delta\psi)\zeta(\psi_{(r)})\varphi+\left\langle \nabla\psi
_{(r)},\nabla\varphi\right\rangle \zeta(\psi_{(r)})]\,d\mu\,dt\nonumber\\
&  =\int[-(\Delta\psi)\zeta(\psi_{(r)})\varphi+(\Delta\varphi)\psi_{(r)}%
\zeta(\psi_{(r)})+\left\langle \nabla\psi,\nabla\varphi\right\rangle
\psi_{(r)}\zeta^{\prime}(\psi_{(r)})]\,d\mu\,dt,\nonumber
\end{align}
the calculation being valid on each time slice.

For the fourth integral, we compute that at each time $t\in(a,b)$, one has%
\[
\frac{d}{dt}\int_{\mathcal{M}^{n}\times\{t\}}Z(\psi_{(r)})\varphi\,d\mu
=\int_{\mathcal{M}^{n}\times\{t\}}[Z^{\prime}(\psi_{(r)})\frac{\partial\psi
}{\partial t}\varphi+Z(\psi_{(r)})\frac{\partial\varphi}{\partial t}%
+Z(\psi_{(r)})\varphi(\operatorname*{tr}\!_{g}h)]\,d\mu,
\]
the final term coming from differentiation of the volume form. Integrating
over the time interval $[a,s]$ and using the facts that $Z^{\prime}=\zeta$ and
that $Z(\psi_{(r)})=0$ at $t=a$ (which holds because $\psi_{(r)}\leq0$ at
$t=a$ by Remark~\ref{CompactHeatballs}), we find that%
\begin{align}
\int Z(\psi_{(r)})(\operatorname*{tr}\!_{g}h)\varphi\,d\mu\,dt  &
=-\int[\zeta(\psi_{(r)})\frac{\partial\psi}{\partial t}\varphi+Z(\psi
_{(r)})\frac{\partial\varphi}{\partial t}]\,d\mu\,dt\label{term4calc}\\
&  +\int_{\mathcal{M}^{n}\times\{s\}}Z(\psi_{(r)})\varphi\,d\mu,\nonumber
\end{align}
where the integrals are still over $\mathcal{M}^{n}\times\lbrack a,s]$ unless
otherwise indicated.

We now combine (\ref{Qderiv1}) with (\ref{term1calc}) and (\ref{term4calc}) to
obtain%
\begin{align}
\frac{r^{n+1}}{n}\frac{d}{dr}\left[  \frac{Q(r)}{r^{n}}\right]   &
=-\int(\frac{\partial\psi}{\partial t}+\Delta\psi+|\nabla\psi|^{2}%
+\operatorname*{tr}\!_{g}h)\zeta(\psi_{(r)})\varphi\,d\mu\,dt\\
&  +\int(\Delta\varphi)\psi_{(r)}\zeta(\psi_{(r)})\,d\mu\,dt-\int
\frac{\partial\varphi}{\partial t}Z(\psi_{(r)})\,d\mu\,dt\nonumber\\
&  +\int\left\langle \nabla\psi,\nabla\varphi\right\rangle \psi_{(r)}%
\zeta^{\prime}(\psi_{(r)})\,d\mu\,dt+\int_{\mathcal{M}^{n}\times\{s\}}%
Z(\psi_{(r)})\varphi\,d\mu.\nonumber
\end{align}
The entire identity may now be multiplied by $n/r^{n+1}$ and integrated with
respect to $r$ between $r_{0}$ and $r_{1}$, where $0<r_{0}<r_{1}\leq\bar{r}$,
to get an identity for the quantity $Q(r_{1},s)/r_{1}^{n}-Q(r_{0},s)/r_{0}%
^{n}$.

We may simplify the resulting expression by picking an appropriate sequence of
valid functions $\zeta$ and passing to the limit. Precisely, we pick a smooth
$\zeta_{1}:\mathbb{R}\rightarrow\lbrack0,1]$ with the properties that
$\zeta_{1}(y)=0$ for $y\leq1/2$, $\zeta_{1}(y)=1$ for $y\geq1$, and $\zeta
_{1}^{\prime}(y)\geq0$. Then we define a sequence $\zeta_{k}:\mathbb{R}%
\rightarrow\lbrack0,1]$ by $\zeta_{k}(y)=\zeta_{1}(2^{k-1}y)$. As $k$
increases, this sequence increases pointwise to the characteristic function of
$(0,\infty)$. The corresponding $Z_{k}$ converge uniformly to the function
$y\mapsto\lbrack y]_{+}$. Crucially, we also can make use of the facts that
$\zeta_{k}(\psi_{(r)})$ converges to the characteristic function of $E_{r}$ in
$L^{1}(\mathcal{M}^{n}\times\lbrack a,b])$ and that $\psi_{(r)}\zeta
_{k}^{\prime}(\psi_{(r)})$ is a bounded sequence of functions on
$\mathcal{M}^{n}\times\lbrack a,b)$ with disjoint supports for each $k$.
Indeed, the support of $\zeta_{k}^{\prime}$ lies within the interval
$(2^{-k},2^{1-k})$.

For each $r\in(0,\bar{r}]$, we have $Q(r,s)\rightarrow P(r,s)$ as
$k\rightarrow\infty$. Using the dominated convergence theorem, our expression
becomes%
\begin{align}
\frac{P(r_{1},s)}{r_{1}^{n}}  &  -\frac{P(r_{0},s)}{r_{0}^{n}}\\
&  =-\int_{r_{0}}^{r_{1}}\frac{n}{r^{n+1}}\int_{E_{r}\cap(\mathcal{M}%
^{n}\times\lbrack a,s])}(\frac{\partial\psi}{\partial t}+\Delta\psi
+|\nabla\psi|^{2}+\operatorname*{tr}\!_{g}h)\varphi\,d\mu\,dt\,dr\nonumber\\
&  -\int_{r_{0}}^{r_{1}}\frac{n}{r^{n+1}}\int_{E_{r}\cap(\mathcal{M}^{n}%
\times\lbrack a,s])}\psi_{(r)}(\frac{\partial\varphi}{\partial t}%
-\Delta\varphi)\,d\mu\,dt\,dr\nonumber\\
&  +\int_{r_{0}}^{r_{1}}\frac{n}{r^{n+1}}\int_{\mathcal{M}^{n}\times
\{s\}}[\psi_{(r)}]_{+}\varphi\,d\mu\,dr.\nonumber
\end{align}
Now we may take the limit as $s\nearrow b$. By Assumption~\ref{r-bar}, the
final term converges to zero, and we end up with (\ref{Pexp1}) as desired.

\medskip

Next we turn to the case that $\psi$ is merely Lipschitz, in the sense of
Assumption~\ref{LL}. Given $s\in(a,b)$ and functions $\zeta$ and $Z$ as above,
there exists a sequence of smooth functions $\psi_{j}$ on $\mathcal{M}%
^{n}\times\lbrack a,b]$ such that $\psi_{j}\rightarrow\psi$ in both $W^{1,2}$
and $C^{0}$ on the set $E_{\bar{r}}\cap(\mathcal{M}^{n}\times\lbrack a,s])$.
By hypothesis on our Lipschitz $\psi$, we have%
\[
I:=\int_{\mathcal{M}^{n}\times\lbrack a,s]}-(\frac{\partial\psi}{\partial
t}+\Delta\psi+|\nabla\psi|^{2}+\operatorname*{tr}\!_{g}h)\zeta(\psi
_{(r)})\varphi\,d\mu\,dt\leq0,
\]
where we make sense of the Laplacian term via integration by parts, namely%
\[
\int-(\Delta\psi)\zeta(\psi_{(r)})\varphi\,d\mu\,dt:=\int[\left\langle
\nabla\psi,\nabla\varphi\right\rangle \zeta(\psi_{(r)})+|\nabla\psi|^{2}%
\zeta^{\prime}(\psi_{(r)})\varphi]\,d\mu\,dt.
\]
By definition of $\psi_{j}$, we have%
\[
\lim_{j\rightarrow\infty}\int_{\mathcal{M}^{n}\times\lbrack a,s]}%
-(\frac{\partial\psi_{j}}{\partial t}+\Delta\psi_{j}+|\nabla\psi_{j}%
|^{2}+\operatorname*{tr}\!_{g}h)\zeta((\psi_{j})_{(r)})\varphi\,d\mu
\,dt=I\leq0,
\]
uniformly in $r\in(0,\bar{r}]$. Consequently, we may carry out the same
calculations that we did in the first part of the proof to obtain an
inequality for the quantity $Q(r_{1},s)/r_{1}^{n}-Q(r_{0},s)/r_{0}^{n}$, with
$\psi_{j}$ in place of $\psi$. We then pass to the limit as $j\rightarrow
\infty$ to obtain the inequality%
\begin{align}
\frac{Q(r_{1},s)}{r_{1}^{n}}-\frac{Q(r_{0},s)}{r_{0}^{n}}  &  \leq\int_{r_{0}%
}^{r_{1}}\frac{n}{r^{n+1}}\int_{\mathcal{M}^{n}\times\lbrack a,s]}%
(\Delta\varphi)\psi_{(r)}\zeta(\psi_{(r)})\,d\mu\,dt\,dr\\
&  -\int_{r_{0}}^{r_{1}}\frac{n}{r^{n+1}}\int_{\mathcal{M}^{n}\times\lbrack
a,s]}\frac{\partial\varphi}{\partial t}Z(\psi_{(r)})\,d\mu\,dt\,dr\nonumber\\
&  +\int_{r_{0}}^{r_{1}}\frac{n}{r^{n+1}}\int_{\mathcal{M}^{n}\times\lbrack
a,s]}\left\langle \nabla\psi,\nabla\varphi\right\rangle \psi_{(r)}%
\zeta^{\prime}(\psi_{(r)})\,d\mu\,dt\,dr\nonumber\\
&  +\int_{r_{0}}^{r_{1}}\frac{n}{r^{n+1}}\int_{\mathcal{M}^{n}\times
\{s\}}Z(\psi_{(r)})\varphi\,d\mu\,dr\nonumber
\end{align}
for our Lipschitz $\psi$. Finally, we replace $\zeta$ with the same sequence
of cut-off functions $\zeta_{k}$ that we used before (thus approximating the
Heaviside function), take the limit as $k\rightarrow\infty$, and then take the
limit as $s\nearrow b$. This gives the inequality (\ref{Pexp2}).
\end{proof}

The argument above may be compared to proofs of earlier results, especially
the proof \cite{Ecker01} of the local monotonicity formula for mean curvature flow.

\bigskip

There is an alternative formula for (\ref{DefineP}) that we find useful in the sequel:

\begin{lemma}
\label{AlternateP}Suppose that $(\mathcal{M}^{n},g(t))$ is a smooth
one-parameter family of complete Riemannian manifolds evolving by
(\ref{GeneralMetricEvolution}) for $t\in\lbrack a,b]$, that $\Psi
:\mathcal{M}^{n}\times\lbrack a,b)\rightarrow(0,\infty)$ satisfies
Assumptions~1--3, that $\bar{r}>0$ is determined by Assumption~\ref{r-bar} and
Remark~\ref{CompactHeatballs}, and that $0<r_{0}<r_{1}\leq\bar{r}$.

If $\varphi\equiv1$ and $\frac{\partial\psi}{\partial t}+|\nabla\psi|^{2}\in
L^{1}(E_{\bar{r}})$, then for all $r\in(0,\bar{r}]$, one has
\[
P_{\varphi,\Psi}(r)=\int_{E_{r}}(\frac{\partial\psi}{\partial t}+|\nabla
\psi|^{2})\,d\mu\,dt.
\]

\end{lemma}

\begin{proof}
In the case that $\varphi\equiv1$, substituting formula (\ref{term4calc}) into
formula (\ref{DefineQ}) yields%
\[
Q(r,s)=\int_{\mathcal{M}^{n}\times\lbrack a,s]}(\frac{\partial\psi}{\partial
t}+|\nabla\psi|^{2})\zeta(\psi_{(r)})\,d\mu\,dt-\int_{\mathcal{M}^{n}%
\times\{s\}}Z(\psi_{(r)})\,d\mu.
\]
Although (\ref{term4calc}) was derived assuming smoothness of $\psi$, one can
verify that it holds for locally Lipschitz $\Psi$ satisfying
Assumption~\ref{LL} by approximating $\psi$ by a sequence of smooth $\psi_{j}$
(as in the proof of Theorem~\ref{MainTheorem}) and then passing to the limit
as $j\rightarrow\infty$. Then if $\frac{\partial\psi}{\partial t}+|\nabla
\psi|^{2}\in L^{1}(E_{\bar{r}})$, one may (again as in the proof of
Theorem~\ref{MainTheorem}) choose a sequence $\zeta_{k}$ along which
$Q(r,s)\rightarrow P(r,s)$ as $k\rightarrow\infty$ and then let $s\nearrow b$
to obtain the stated formula.
\end{proof}

\section{A local gradient estimate\label{LGE}}

In order to apply Theorem~\ref{MainTheorem} to a fundamental solution of the
heat equation of an evolving manifold in Section \ref{MVT}, we need a local
gradient estimate. One approach would be to adapt existing theory of local
heat kernel asymptotics. Instead, we prove a more general result which may be
of independent interest. Compare \cite{Ham93}, \cite{LY}, \cite{MS96}, the
recent \cite{SZ}, \cite{ST98}, and \cite{Yau75, Yau95}.

\medskip

Let $(\mathcal{M}^{n},g(t))$ be a smooth one-parameter family of complete
Riemannian manifolds evolving by (\ref{GeneralMetricEvolution}) for
$t\in\lbrack0,\bar{t}]$. We shall abuse notation by writing $g(\tau)$ to mean
$g(\tau(t))$, where
\[
\tau(t):=\bar{t}-t.
\]
In the remainder of this section, we state our results solely in terms of
$\tau$. In particular, $g(\tau)$ satisfies $\frac{\partial}{\partial\tau
}g=-2h$ on $\mathcal{M}^{n}\times\lbrack0,\bar{t}]$.

Given $\bar{x}\in\mathcal{M}^{n}$ and $\rho>0$, define%
\begin{equation}
\Omega(\rho):=%
{\displaystyle\bigcup\limits_{0\leq\tau\leq\bar{t}}}
\left(  B_{g(\tau)}(\bar{x},\rho)\times\{\tau\}\right)  \subseteq
\mathcal{M}^{n}\times\lbrack0,\bar{t}]. \label{Omega-rho}%
\end{equation}
We now prove a local \emph{a priori }estimate for bounded positive solutions
of the conjugate heat equation
\begin{equation}
(\frac{\partial}{\partial\tau}-\Delta-\operatorname*{tr}\!_{g}h)v=0.
\label{BCHE}%
\end{equation}
We will apply this in Section~\ref{MVT}.

\begin{theorem}
\label{LocalGradientEstimate}Let $(\mathcal{M}^{n},g(\tau))$ be a smooth
one-parameter family of complete Riemannian manifolds evolving by
$\frac{\partial}{\partial\tau}g=-2h$ for $0\leq\tau\leq\bar{t}$. Assume there
exist $k_{1},k_{2},k_{3}\geq0$ such that
\[
h\leq k_{1}g,\qquad\operatorname*{Rc}\geq-k_{2}g,\qquad\text{and}\qquad
|\nabla(\operatorname*{tr}\!_{g}h)|\leq k_{3}%
\]
in the space-time region $\Omega(2\rho)$ given by (\ref{Omega-rho}). Assume
further that $v(\tau)$ solves (\ref{BCHE}) and satisfies $0<v\leq A$ in
$\Omega(2\rho)$.

Then there exist a constant $C_{1}$ depending only on $n$ and an absolute
constant $C_{2}$ such that at all $(x,\tau)\in\Omega(\rho)$, one has%
\[
\frac{|\nabla v|^{2}}{v^{2}}\leq(1+\log\frac{A}{v})^{2}\left[  \frac{1}{\tau
}+C_{1}k_{1}+2k_{2}+k_{3}+\sqrt{k_{3}}+\frac{C_{1}\sqrt{k_{2}}\rho\coth
(\sqrt{k_{2}}\rho)+C_{2}}{\rho^{2}}\right]  .
\]

\end{theorem}

\begin{proof}
By scaling, we may assume that $A=1$. We define\thinspace\footnote{Note that
$w$ is used by Souplet--Zhang \cite[Theorem 1.1]{SZ} in generalizing
Hamilton's result \cite{Ham93}. A similar function is employed by Yau
\cite{Yau75}. Also see related work of the third author \cite{Ni-PLYH}.}%
\[
f:=\log v\qquad\text{and}\qquad w:=|\nabla\log(1-f)|^{2},
\]
computing that
\[
(\frac{\partial}{\partial\tau}-\Delta)f=|\nabla f|^{2}+(\operatorname*{tr}%
\!_{g}h).
\]
Then using Bochner--Weitzenb\"{o}ck, we calculate that%
\[
(\frac{\partial}{\partial\tau}-\Delta)|\nabla f|^{2}=2h(\nabla f,\nabla
f)-2\operatorname*{Rc}(\nabla f,\nabla f)-2|\nabla\nabla f|^{2}+2\left\langle
\nabla(\operatorname*{tr}\!_{g}h+|\nabla f|^{2}),\nabla f\right\rangle
\]
and%
\begin{align*}
(\frac{\partial}{\partial\tau}-\Delta)w  &  =\frac{1}{(1-f)^{2}}[2h(\nabla
f,\nabla f)-2\operatorname*{Rc}(\nabla f,\nabla f)+2\left\langle
\nabla(\operatorname*{tr}\!_{g}h),\nabla f\right\rangle ]\\
&  -\frac{2}{(1-f)^{2}}\left[  |\nabla\nabla f|^{2}+\frac{\left\langle
\nabla|\nabla f|^{2},\nabla f\right\rangle }{1-f}+\frac{|\nabla f|^{4}%
}{(1-f)^{2}}\right] \\
&  +2\frac{(\operatorname*{tr}\!_{g}h)|\nabla f|^{2}+|\nabla f|^{4}}%
{(1-f)^{3}}-4\frac{|\nabla f|^{4}}{(1-f)^{4}}-2f\frac{\left\langle
\nabla|\nabla f|^{2},\nabla f\right\rangle }{(1-f)^{3}}.
\end{align*}
Rewriting some factors $\left\langle \nabla|\nabla f|^{2},\nabla
f\right\rangle $ in terms of $\left\langle \nabla w,\nabla f\right\rangle $,
we obtain%
\begin{align*}
(\frac{\partial}{\partial\tau}-\Delta)w  &  =\frac{1}{(1-f)^{2}}[2h(\nabla
f,\nabla f)-2\operatorname*{Rc}(\nabla f,\nabla f)+2\left\langle
\nabla(\operatorname*{tr}\!_{g}h),\nabla f\right\rangle ]\\
&  -\frac{2}{(1-f)^{2}}\left[  |\nabla\nabla f|^{2}+\frac{\left\langle
\nabla|\nabla f|^{2},\nabla f\right\rangle }{1-f}+\frac{|\nabla f|^{4}%
}{(1-f)^{2}}\right] \\
&  +2\frac{(\operatorname*{tr}\!_{g}h)|\nabla f|^{2}-|\nabla f|^{4}}%
{(1-f)^{3}}+2\frac{-f}{1-f}\left\langle \nabla w,\nabla f\right\rangle .
\end{align*}

Now let $\eta(s)$ be a smooth nonnegative cutoff function such that
$\eta(s)=1$ when $s\leq1$ and $\eta(s)=0$ when $s\geq2$, with $\eta^{\prime
}\leq0$, $|\eta^{\prime}|\leq C_{2}$, $(\eta^{\prime})^{2}\leq C_{2}\eta$, and
$\eta^{\prime\prime}\geq-C_{2}$. Define%
\[
u(x,\tau):=\eta\left(  \frac{d_{g(\tau)}(\bar{x},x)}{\rho}\right)  .
\]
Then we have%
\[
\frac{|\nabla u|^{2}}{u}\leq\frac{C_{2}}{\rho^{2}}%
\]
and%
\[
\frac{\partial u}{\partial\tau}\leq C_{2}k_{1}%
\]
and%
\[
-\Delta u\leq\frac{C_{1}\sqrt{k_{2}}\rho\coth(\sqrt{k_{2}}\rho)+C_{2}}%
{\rho^{2}}.
\]
These inequalities hold in the barrier sense. However, when applying the
maximum principle, Calabi's trick lets us pretend that $u$ is actually smooth.

Now let $G:=uw$ and compute that%
\[
(\frac{\partial}{\partial\tau}-\Delta)(\tau G)=G+\tau u\left[  (\frac
{\partial}{\partial\tau}-\Delta)w\right]  +\tau w\left[  (\frac{\partial
}{\partial\tau}-\Delta)u\right]  -2\tau\left\langle \nabla u,\nabla
w\right\rangle .
\]
For any $\tau_{1}>0$, consider $\tau G$ on $\mathcal{M}^{n}\times\lbrack
0,\tau_{1}]$. At any point $(x_{0},\tau_{0})$ where $\tau G$ attains its
maximum on $\mathcal{M}^{n}\times\lbrack0,\tau_{1}]$, we have $0\leq
(\frac{\partial}{\partial\tau}-\Delta)(\tau G)$ and%
\begin{align*}
(\frac{\partial}{\partial\tau}-\Delta)(\tau G)  &  \leq G-2\tau\left\langle
\nabla u,\nabla w\right\rangle \\
&  +2\tau u\left[  \frac{h(\nabla f,\nabla f)-\operatorname*{Rc}(\nabla
f,\nabla f)+\left\langle \nabla(\operatorname*{tr}\!_{g}h),\nabla
f\right\rangle }{(1-f)^{2}}+\frac{(\operatorname*{tr}\!_{g}h)|\nabla f|^{2}%
}{(1-f)^{3}}\right] \\
&  +2\tau u\left[  \frac{-f}{1-f}\left\langle \nabla w,\nabla f\right\rangle
-\frac{|\nabla f|^{4}}{(1-f)^{3}}\right] \\
&  +\tau w\left[  C_{2}k_{1}+\frac{C_{1}\sqrt{k_{2}}\rho\coth(\sqrt{k_{2}}%
\rho)+C_{2}}{\rho^{2}}\right]  .
\end{align*}
Using the fact that $\nabla G(x_{0},t_{0})=0$, we can replace $u\nabla w$ by
$-w\nabla u$ above. Then multiplying both sides of the inequality by
$u\in\lbrack0,1]$ and using $1/(1-f)\leq1$, we obtain%
\begin{align*}
0  &  \leq G+2\tau\left\{  \lbrack(n+1)k_{1}+k_{2}]G+k_{3}\sqrt{G}\right\} \\
&  +2\tau G|\nabla u||\nabla f|\left(  \frac{-f}{1-f}\right)  -2\tau
(1-f)G^{2}\\
&  +\tau G\left\{  C_{2}k_{1}+\frac{C_{1}\sqrt{k_{2}}\rho\coth(\sqrt{k_{2}%
}\rho)+C_{2}}{\rho^{2}}\right\}  .
\end{align*}
Noticing that $2k_{3}\sqrt{G}\leq k_{3}(G+1)$ and that%
\begin{align*}
2\tau G|\nabla u||\nabla f|\left(  \frac{-f}{1-f}\right)   &  \leq\tau
G\left(  \frac{|\nabla f|^{2}}{1-f}u+\frac{|\nabla u|^{2}}{u}\frac{f^{2}}%
{1-f}\right) \\
&  \leq\tau(1-f)G^{2}+\tau G\frac{C_{2}}{\rho^{2}}\frac{f^{2}}{1-f},
\end{align*}
we estimate at $(x_{0},t_{0})$ that%
\begin{align*}
0  &  \leq\tau k_{3}+G\left\{  1+\tau\left[  C_{1}k_{1}+2k_{2}+k_{3}%
+\frac{C_{1}\sqrt{k_{2}}\rho\coth(\sqrt{k_{2}}\rho)+C_{2}}{\rho^{2}}\right]
\right\} \\
&  +\tau G\frac{C_{2}}{\rho^{2}}\frac{f^{2}}{1-f}-\tau(1-f)G^{2}.
\end{align*}
Dividing both sides by $\tau(1-f)$ while noting that $1/(1-f)\leq1$ and
$-f/(1-f)\leq1$, we get%
\[
0\leq k_{3}+G\left[  \frac{1}{\tau}+C_{1}k_{1}+2k_{2}+k_{3}+\frac{C_{1}%
\sqrt{k_{2}}\rho\coth(\sqrt{k_{2}}\rho)+C_{2}}{\rho^{2}}\right]  -G^{2},
\]
from which we can conclude that%
\[
G\leq\frac{1}{\tau}+C_{1}k_{1}+2k_{2}+k_{3}+\sqrt{k_{3}}+\frac{C_{1}%
\sqrt{k_{2}}\rho\coth(\sqrt{k_{2}}\rho)+C_{2}}{\rho^{2}}%
\]
at $(x_{0},\tau_{0})$. Hence $W(\tau_{1}):=\tau_{1}\sup_{x\in B_{g(\tau)}%
(\bar{x},\rho)}w(x,\tau_{1})$ may be estimated by%
\begin{align*}
W(\tau_{1})  &  \leq\tau_{0}G(x_{0},\tau_{0})\\
&  \leq1+\tau_{0}\left[  C_{1}k_{1}+2k_{2}+k_{3}+\sqrt{k_{3}}+\frac{C_{1}%
\sqrt{k_{2}}\rho\coth(\sqrt{k_{2}}\rho)+C_{2}}{\rho^{2}}\right] \\
&  \leq1+\tau_{1}\left[  C_{1}k_{1}+2k_{2}+k_{3}+\sqrt{k_{3}}+\frac{C_{1}%
\sqrt{k_{2}}\rho\coth(\sqrt{k_{2}}\rho)+C_{2}}{\rho^{2}}\right]  .
\end{align*}
Since $\tau_{1}>0$ was arbitrary, the result follows.
\end{proof}

\begin{remark}
In the special case that $h\equiv0$, we have%
\[
\frac{|\nabla v|^{2}}{v^{2}}\leq(1+\log\frac{A}{v})^{2}\left[  \frac{1}{\tau
}+2k_{2}+\frac{C_{1}\sqrt{k_{2}}\rho\coth(\sqrt{k_{2}}\rho)+C_{2}}{\rho^{2}%
}\right]
\]
at $(x,\tau)$, for all times $\tau\in\lbrack0,\bar{t}]$ and points $x\in
B_{g(\tau)}(\bar{x},\rho)$, which slightly improves a result of \cite{SZ}.
\end{remark}

\section{Reduced volume for Ricci flow\label{RV}}

Our first application of Theorem~\ref{MainTheorem} is to Ricci flow. Let
$(\mathcal{M}^{n},g(t))$ be a complete solution of Ricci flow that remains
smooth for $0\leq t\leq\bar{t}$. This corresponds to $h=-\operatorname*{Rc}$
and $\operatorname*{tr}\!_{g}h=-R$ in (\ref{GeneralMetricEvolution}).

\subsection{Localizing Perelman's reduced volume}

Perelman \cite[\S 7]{Perelman1} has discovered a remarkable quantity that may
be regarded as a kind of parabolic distance for Ricci flow. Define
$\tau(t):=\bar{t}-t$, noting that $g(\tau(t))$ then satisfies $\frac{\partial
}{\partial\tau}g=2\operatorname*{Rc}$ for $0\leq\tau\leq\bar{t}$. Fix $\bar
{x}\in\mathcal{M}^{n}$ and regard $(\bar{x},0)$ (in $(x,\tau)$ coordinates) as
a space-time origin. The space-time action of a smooth path $\gamma$ with
$\gamma(0)=(\bar{x},0)$ and $\gamma(\tau)=(x,\tau)$ is
\begin{subequations}
\label{ScriptL}%
\begin{align}
\mathcal{L}(\gamma)  &  :=\int_{0}^{\tau}\sqrt{\sigma}\left(  |\frac{d\gamma
}{d\sigma}|^{2}+R\right)  \,d\sigma\\
&  =\int_{0}^{\sqrt{\tau}}\left(  \frac{1}{2}|\frac{d\gamma}{ds}|^{2}%
+2s^{2}R\right)  \,ds\qquad\qquad(s=\sqrt{\sigma}).
\end{align}
Taking the infimum over all such paths, Perelman defines the reduced distance
from $(\bar{x},0)$ to $(x,\tau)$ as%
\end{subequations}
\begin{equation}
\ell(x,\tau)=\ell_{(\bar{x},0)}(x,\tau):=\frac{1}{2\sqrt{\tau}}\inf_{\gamma
}\mathcal{L}(\gamma), \label{ReducedDistance}%
\end{equation}
and observes that
\begin{equation}
v(x,\tau)=v_{(\bar{x},0)}(x,\tau):=\frac{1}{(4\pi\tau)^{n/2}}e^{-\ell(x,\tau)}
\label{ReducedVolumeDensity}%
\end{equation}
is a subsolution of the conjugate heat equation $u_{\tau}=\Delta u-Ru$ in the
barrier sense \cite{Perelman1}, hence in the distributional
sense.\footnote{See \cite{Ye} for a direct proof of the distributional
inequality.} It follows that the reduced volume (essentially a Gaussian
weighted volume)
\begin{equation}
\tilde{V}(t)=\tilde{V}_{(\bar{x},0)}(t):=\int_{\mathcal{M}^{n}}v(x,\tau)\,d\mu
\end{equation}
is a monotonically increasing function of $t$ which is constant precisely on
shrinking gradient solitons. (Compare to Example~\ref{RVE} above.)

The interpretations of $\ell$ as parabolic distance and $\tilde{V}$ as
Gaussian weighted volume are elucidated by the following examples.

\begin{example}
Let $(\mathcal{M}^{n},g)$ be a Riemannian manifold of nonnegative Ricci
curvature, and let $q$ be any smooth superharmonic function $(\Delta q\leq0)$.
In their seminal paper \cite{LY}, Li and Yau define%
\[
\rho(x,\tau)=\inf\left\{  \frac{1}{4\tau}\int_{0}^{1}|\frac{d\gamma}{d\sigma
}|^{2}\,d\sigma+\tau\int_{0}^{1}q(\gamma(\sigma))\,d\sigma\right\}  ,
\]
where the infimum is taken over all smooth paths from an origin $(\bar{x},0)$.
As a special case of their more general results \cite[Theorem~4.3]{LY}, they
observe that $(4\pi\tau)^{-n/2}e^{-\rho(x,\tau)}$ is a distributional
subsolution of the linear parabolic equation $(\frac{\partial}{\partial\tau
}-\Delta+q)u=0$.
\end{example}

\begin{example}
\label{RV1}Let $(\mathbb{R}^{n},g)$ denote Euclidean space with its standard
flat metric. Given $\lambda\in\mathbb{R}$, define $X=\operatorname*{grad}%
(\frac{\lambda}{2}|x|^{2}|)$. Then one has $0=\operatorname*{Rc}=\lambda
g-\mathcal{L}_{X}g$. Hence there is a Ricci soliton structure on Euclidean
space, called the Gaussian soliton.

Take $\lambda=1$ to give $(\mathbb{R}^{n},g)$ the structure of a gradient
shrinking soliton. Then $\gamma(\sigma)=\sqrt{\sigma/\tau}\,x$ is an
$\mathcal{L}$-geodesic from $(0,0)$ to $(x,\tau)$. Thus the reduced distance
is $\ell_{(0,0)}(x,\tau)=|x|^{2}/4\tau$ and the reduced volume integrand is
exactly the heat kernel $v_{(0,0)}(x,\tau)=(4\pi\tau)^{-n/2}e^{-|x|^{2}/4\tau
}$. Hence $\tilde{V}_{(0,0)}(t)\equiv1$. (Compare \cite[\S 15]{KL}.)
\end{example}

\begin{example}
\label{RV2}Let $S_{r(\tau)}^{n}$ denote the round sphere of radius
$r(\tau)=\sqrt{2(n-1)\tau}$. This is a positive Einstein manifold, hence a
homothetically shrinking (in $t$) solution of Ricci flow. Along any sequence
$(x_{k},\tau_{k})$ of smooth origins approaching the singularity $\mathcal{O}$
at $\tau=0$, one gets a smooth function $\ell_{\mathcal{O}}(x,\tau
):=\lim_{k\rightarrow\infty}\ell_{(x_{k},\tau_{k})}(x,\tau)\equiv n/2$
measuring the reduced distance from $\mathcal{O}$. Hence $\tilde
{V}_{\mathcal{O}}(t)\equiv\lbrack(n-1)/(2\pi e)]^{n/2}\operatorname*{Vol}%
(S_{1}^{n})$ for all $t<0$.
\end{example}

Our first application of Theorem~\ref{MainTheorem} is where $\Psi$ is
Perelman's reduced-volume density $v$ (\ref{ReducedVolumeDensity}). Let $\ell$
denote the reduced distance (\ref{ReducedDistance}) from a smooth origin
$(\bar{x},\bar{t})$ and assume there exists $k\in(0,\infty)$ such that
$\operatorname*{Rc}\geq-kg$ on $\mathcal{M}^{n}\times\lbrack0,\bar{t}]$. In
what follows, we will freely use results from the Appendix
(Section~\ref{Appendix}, below).

Lemma~\ref{Lipschitz} guarantees that $\ell$ is locally Lipschitz, hence that
Assumption~\ref{LL} is satisfied. (Also see \cite{Ye} or \cite{RFV2}.) The
estimate in Part~(1) of Lemma~\ref{BoundsForReducedDistance} ensures that
Assumption~\ref{Omega} is satisfied. Assumption~\ref{r-bar} follows from
combining that estimate, Corollary~\ref{BoundReducedBalls}, and Lemma
\ref{Finiteness}. Here we may take any $\bar{r}>0$ satisfying $\bar{r}%
^{2}<\min\{\bar{t}/c,4\pi\}$, where $c=e^{4k\bar{t}/3}/(4\pi)$. So for
$r\in(0,\bar{r}]$, consider
\[
P_{\varphi,v}(r):=\int_{E_{r}}[|\nabla\ell|^{2}+R(n\log\frac{r}{\sqrt{4\pi
\tau}}-\ell)]\varphi\,d\mu\,dt.
\]
Notice that $|\nabla\ell|^{2}$ replaces the term $\frac{|x-\bar{x}|^{2}}%
{4\tau^{2}}$ in the heatball formulas for Euclidean space and solutions of
mean curvature flow. See Examples~\ref{Euclidean} and \ref{MCFdensity}, respectively.

\begin{remark}
\label{RD-alt}For $r\in(0,\bar{r}]$, one may write $P_{1,v}(r)$ in either
alternative form
\begin{subequations}
\label{AlternateReduced}%
\begin{align}
P_{1,v}(r)  &  =\int_{E_{r}}(\frac{n}{2\tau}+\ell_{\tau}+|\nabla\ell
|^{2})\,d\mu\,dt\\
&  =\int_{E_{r}}(\frac{n}{2\tau}-\frac{1}{2}\tau^{-3/2}\mathcal{K})\,d\mu\,dt.
\end{align}
Here $\mathcal{K}(x,\tau)=\int_{0}^{\tau}\sigma^{3/2}\mathcal{H}%
(d\gamma/d\sigma)\,d\sigma$ is computed along a minimizing $\mathcal{L}%
$-geodesic $\gamma$, where $\mathcal{H}(X)=2\operatorname*{Rc}(X,X)-(R_{\tau
}+2\left\langle \nabla R,X\right\rangle +R/\tau)$ is Hamilton's traced
differential Harnack expression.

If $R\geq0$ and $\varphi\geq0$ on $E_{\bar{r}}$, then for all $r\in(0,\bar
{r}]$, one has%
\end{subequations}
\begin{equation}
P_{\varphi,v}(r)=\int_{E_{r}}[|\nabla\ell|^{2}+R\psi_{(r)}]\varphi
\,d\mu\,dt\geq\int_{E_{r}}|\nabla\ell|^{2}\varphi\,d\mu\,dt\geq0.
\label{PositiveP}%
\end{equation}

If $(\mathcal{M}^{n},g(0))$ has nonnegative curvature operator and $\bar
{r}^{2}<4\pi\bar{t}(1-1/C)$ for some $C>1$, then for all $r\in(0,\bar{r}]$,%
\begin{equation}
P_{1,v}(r)\leq\int_{E_{r}}\frac{n/2+C\ell}{\tau}\,d\mu\,dt. \label{EstimateP}%
\end{equation}

\end{remark}

\begin{proof}
By Part~(2) of Lemma~\ref{Lipschitz}, the arguments of Lemma~\ref{Finiteness}
apply to show that $\psi_{t}+|\nabla\psi|^{2}=\frac{n}{2\tau}+\ell_{\tau
}+|\nabla\ell|^{2}\in L^{1}(E_{\bar{r}})$. Hence Lemma~\ref{AlternateP} and
identities~(7.5) and (7.6) of \cite{Perelman1} imply
formulae~(\ref{AlternateReduced}).

Since $\psi_{(r)}>0$ in $E_{r}$, the inequalities in (\ref{PositiveP}) are clear.

If $(\mathcal{M}^{n},g(0))$ has nonnegative curvature operator, Hamilton's
traced differential Harnack inequality \cite{Ham93b} implies that
\[
\mathcal{H}(\frac{d\gamma}{d\sigma})\geq-R(\frac{1}{\sigma}+\frac{1}{\bar
{t}-\sigma})=-\frac{\bar{t}}{\bar{t}-\sigma}\frac{R}{\sigma}%
\]
along a minimizing $\mathcal{L}$-geodesic $\gamma$. Hence%
\[
-\frac{1}{2}\tau^{-3/2}\mathcal{K\leq}\frac{\bar{t}}{\bar{t}-\tau}\frac
{\tau^{-3/2}}{2}\int_{0}^{\tau}\sqrt{\sigma}(R+|\frac{d\gamma}{d\sigma}%
|^{2})\,d\sigma=\frac{\bar{t}}{\bar{t}-\tau}\frac{\ell}{\tau}.
\]
By Lemma~\ref{rhoEstimates}, one has $\tau<r^{2}/4\pi$, which gives
estimate~(\ref{EstimateP}).
\end{proof}

\bigskip

Our main result in this section is as follows. Recall that $\psi_{(r)}%
:=n\log(\frac{r}{\sqrt{4\pi\tau}})-\ell$.

\begin{corollary}
\label{ReducedDistanceCorollary}Let $(\mathcal{M}^{n},g(t))$ be a complete
solution of Ricci flow that remains smooth for $0\leq t\leq\bar{t}$ and
satisfies $\operatorname*{Rc}\geq-kg$. Let $\varphi$ be any smooth nonnegative
function of $(x,t)$ and let $c=e^{4k\bar{t}/3}/(4\pi)$. Then whenever
$0<r_{0}<r_{1}<\min\{\sqrt{\bar{t}/c},2\sqrt{\pi}\}$, one has%
\begin{equation}
\frac{P_{\varphi,v}(r_{1})}{r_{1}^{n}}-\frac{P_{\varphi,v}(r_{0})}{r_{0}^{n}%
}\leq-\int_{r_{0}}^{r_{1}}\frac{n}{r^{n+1}}\int_{E_{r}}\psi_{(r)}%
(\frac{\partial\varphi}{\partial t}-\Delta\varphi)\,d\mu\,dt\,dr.
\label{ReducedVolumeComparison}%
\end{equation}
Furthermore,%
\begin{equation}
\varphi(\bar{x},\bar{t})=\lim_{r\searrow0}\frac{P_{\varphi,v}(r)}{r^{n}}.
\end{equation}
In particular,%
\begin{equation}
\varphi(\bar{x},\bar{t})\geq\frac{P_{\varphi,v}(r_{1})}{r_{1}^{n}}+\int
_{0}^{r_{1}}\frac{n}{r^{n+1}}\int_{E_{r}}\psi_{(r)}[(\frac{\partial}{\partial
t}-\Delta)\varphi]\,d\mu\,dt\,dr. \label{Klaus}%
\end{equation}

\end{corollary}

\begin{proof}
The quantity $\Psi=$ $v$ satisfies $\frac{\partial v}{\partial t}+\Delta
v-Rv\geq0$ as a distribution. (This is implied by Perelman's barrier
inequality \cite[(7.13)]{Perelman1}; see \cite[Lemma 1.12]{Ye} for a direct
proof.) Hence we may apply Theorem~\ref{MainTheorem} in the form (\ref{Pexp2})
to obtain~(\ref{ReducedVolumeComparison}).

Formula (7.6) of Perelman \cite{Perelman1} implies that%
\[
P_{\varphi,v}(r)=\int_{E_{r}}[\frac{\ell}{\tau}+R\psi_{(r)}-R-\tau
^{-3/2}\mathcal{K}]\varphi\,d\mu\,dt.
\]
By Corollary~\ref{BoundReducedBalls}, there is a precompact neighborhood
$\mathcal{U}$ of $\bar{x}$ with $E_{r}\subseteq\mathcal{U}\times
\lbrack0,cr^{2}]$ for all $r>0$ under consideration. By
Lemma~\ref{Neighborhood}, there exists a precompact set $\mathcal{V}$ such
that the images of all minimizing $\mathcal{L}$-geodesics from $(\bar{x},0)$
to points in $\mathcal{U}\times\lbrack0,cr^{2}]$ are contained in the set
$\mathcal{V}\times\lbrack0,cr^{2}]$, in which one has uniform bounds on all
curvatures and their derivatives. So by Lemma~\ref{BoundsForReducedDistance},
one has $\frac{\ell}{\tau}=\frac{d_{0}^{2}(\bar{x},x)}{4\tau^{2}}+O(\frac
{1}{\tau})$ and $R\psi_{(r)}=R(n\log\frac{r}{\sqrt{4\pi\tau}}-\ell)=O(\frac
{1}{\tau})$ as $\tau\searrow0$. By Corollary~\ref{GeodesicSpeedImproved},
$\tau^{-3/2}\mathcal{K}$ is also $O(\frac{1}{\tau})$ as $\tau\searrow0$.
Adapting the arguments in the proof of Lemma~\ref{Finiteness}, one concludes
that%
\[
\lim_{r\searrow0}\frac{P_{\varphi,v}(r)}{r^{n}}=\lim_{r\searrow0}\left\{
\frac{1}{r^{n}}\int_{E_{r}}\frac{d_{0}^{2}(\bar{x},x)}{4\tau^{2}}\varphi
\,d\mu\,dt\right\}  =\varphi(\bar{x},\bar{t}),
\]
exactly as in the calculation for Euclidean space. (Also see
Corollary~\ref{Moving}, below.)
\end{proof}

An example of how this result may be applied is the following local Harnack
inequality, which follows directly from~(\ref{Klaus}).

\begin{remark}
Assume the hypotheses of Corollary~\ref{ReducedDistanceCorollary} hold. If
$R\geq0$ on $E_{r_{1}}$, then%
\[
R(\bar{x},\bar{t})\geq\frac{1}{r_{1}^{n}}\int_{E_{r_{1}}}[|\nabla\ell
|^{2}+R\psi_{(r)}]R\,d\mu\,dt+\int_{0}^{r_{1}}\frac{2n}{r^{n+1}}\int_{E_{r}%
}\psi_{(r)}\left\vert \operatorname*{Rc}\right\vert ^{2}\,d\mu\,dt\,dr.
\]

\end{remark}

The inequality (\ref{ReducedVolumeComparison}) is sharp in the following sense.

\begin{corollary}
\label{VolumeEquality}Let $(\mathcal{M}^{n},g(t))$ be a complete solution of
Ricci flow that is smooth for $0\leq t\leq\bar{t}$, with $\operatorname*{Rc}%
\geq-kg$. If equality holds in (\ref{ReducedVolumeComparison}) for
$\varphi\equiv1$, then $(E_{r},g(t))$ is isometric to a shrinking gradient
soliton for all $r<\min\{\sqrt{\bar{t}/c},2\sqrt{\pi}\}$.
\end{corollary}

\begin{proof}
From the proof of Theorem~\ref{MainTheorem}, it is easy to see that%
\[
\frac{d}{dr}\left(  \frac{P_{1,v}(r)}{r^{n}}\right)  =-\frac{n}{r^{n+1}}%
\int_{E_{r}}\frac{(\frac{\partial}{\partial t}+\Delta+\operatorname*{tr}%
\!_{g}h)v}{v}\,d\mu\,dt
\]
for almost all $r<\min\{\sqrt{\bar{t}/c},2\sqrt{\pi}\}$. Therefore, equality
in (\ref{ReducedVolumeComparison}) implies that $v$ is a distributional
solution of the parabolic equation%
\[
(\frac{\partial}{\partial\tau}-\Delta+R)v=0
\]
in $E_{r}$ for almost all small $r$. By parabolic regularity, $v$ is actually
smooth. This implies that one has equality in the chain of inequalities%
\[
\Delta\ell-|\nabla\ell|^{2}+R-\frac{n}{2\tau}\leq\ell_{\tau}=-(-R+\frac
{n}{2\tau}-\frac{1}{2}\tau^{-3/2}\mathcal{K})+\frac{n-2\ell}{2\tau}\leq
-\Delta\ell+\frac{n-2\ell}{2\tau}%
\]
that follow from equations (7.13), (7.5), and (7.10) of \cite{Perelman1}.
Hence one has%
\[
u:=\tau(2\Delta\ell-|\nabla\ell|^{2}+R)+\ell-n=0.
\]
By equation (9.1) of \cite{Perelman1} (where the roles of $u$ and $v$ are
reversed), this implies that
\[
0=(\frac{\partial}{\partial\tau}-\Delta+R)(uv)=-2\tau|\operatorname*{Rc}%
+\nabla\nabla\ell-\frac{1}{2\tau}g|^{2}v.
\]
This is possible only if $(E_{r},g(t))$ has the structure of a shrinking
gradient soliton with potential function $\ell$.
\end{proof}

\begin{remark}
For applications of Corollary~\ref{VolumeEquality} to regularity theorems for
Ricci flow, see \cite{Ni-MVT} by the third author.
\end{remark}

\subsection{Comparing global and local quantities}

Corollaries~\ref{ReducedDistanceCorollary} and \ref{VolumeEquality} suggest a
natural question: how does the purely local monotone quantity $P_{1,v}%
(r)/r^{n}$ compare to Perelman's global monotone quantity $\tilde{V}%
(t)=\int_{\mathcal{M}^{n}}v\,d\mu$? A path to a partial answer begins with an
observation that generalizes Example~\ref{RV2} above.

Cao, Hamilton, and Ilmanen \cite{CHI} prove that any complete gradient
shrinking soliton $(\mathcal{M}^{n},g(t))$ that exists up to a maximal time
$T<\infty$ and satisfies\ certain noncollapsing and curvature decay hypotheses
converges as $t\nearrow T$ to an incomplete (possibly empty) metric cone
$(\mathcal{C},d)$, which is smooth except at the parabolic vertex
$\mathcal{O}$. The convergence is smooth except on a compact set (possibly all
of $\mathcal{M}^{n}$) that vanishes into the vertex.\thinspace\footnote{See
\cite{FIK} for examples where $(\mathcal{C},d)=\lim_{\tau\searrow
0}(\mathcal{M}^{n},g(\tau))$ is nonempty.} Furthermore, they prove that along
a sequence $(x_{k},\tau_{k})$ approaching $\mathcal{O}$, a limit
$\ell_{\mathcal{O}}(x,\tau):=\lim\ell_{(x_{k},\tau_{k})}(x,\tau)$ exists for
all $x\in\mathcal{M}^{n}$ and $\tau(t)>0$. They show that the central density
function
\[
\Theta_{\mathcal{O}}^{\mathrm{RF}}(t):=\tilde{V}_{\mathcal{O}}(t)=\lim
_{k\rightarrow\infty}\tilde{V}_{(x_{k},\tau_{k})}(t)
\]
of the parabolic vertex $\mathcal{O}$ is independent of time and satisfies
$\Theta_{\mathcal{O}}^{\mathrm{RF}}(t)\equiv e^{\nu}$, where $\nu$ is the
constant entropy of the soliton $(\mathcal{M}^{n},g(\tau))$.

On a compact soliton, there is a pointwise version of the
Cao--Hamilton--Ilmanen result, due to Bennett Chow and the third author:

\begin{lemma}
\label{CLN-Observation}If $(\mathcal{M}^{n},g(\tau))$ is a compact shrinking
(necessarily gradient) soliton, then the limit $\ell_{\mathcal{O}}(x,\tau)$
exists for all $x\in\mathcal{M}^{n}$ and $\tau(t)>0$. This limit agrees up to
a constant with the soliton potential function $f(x,\tau)$.
\end{lemma}

See \cite{RFV2} for a proof.

\smallskip

Recall that the entropy of a compact Riemannian manifold $(\mathcal{M}^{n},g)$
is%
\[
\nu(\mathcal{M}^{n},g):=\inf\left\{  \mathcal{W}\left(  g,f,\tau\right)  :f\in
C_{0}^{\infty},~\tau>0,~\int_{\mathcal{M}^{n}}(4\pi\tau)^{-n/2}e^{-f}%
\,d\mu=1\right\}  ,
\]
where%
\begin{equation}
\mathcal{W}\left(  g,f,\tau\right)  :\mathcal{=}\int_{\mathcal{M}^{n}}\left[
\tau(|\nabla f|^{2}+R)+f-n\right]  (4\pi\tau)^{-n/2}e^{-f}\,d\mu.
\label{EntropyFormula}%
\end{equation}
(Compare to Example~\ref{Entropy}.) Under the coupled system
\begin{subequations}
\label{RSS}%
\begin{align}
\frac{\partial}{\partial t}g  &  =-2\operatorname*{Rc}\\
(\frac{\partial}{\partial t}+\Delta)f  &  =|\nabla f|^{2}-R+\frac{n}{2\tau}\\
\frac{d\tau}{dt}  &  =-1,
\end{align}
the functional $\mathcal{W}(g(t),f(t),\tau(t))$ is monotone increasing in time
and is constant precisely on a compact shrinking gradient soliton with
potential function $f$, where (after possible normalization) one has%
\end{subequations}
\begin{equation}
\operatorname*{Rc}+\nabla\nabla f-\frac{1}{2\tau}g\circeq0. \label{RSE}%
\end{equation}
Here and in the remainder of this section, the symbol $\circeq$ denotes an
identity that holds on a shrinking gradient soliton.

We are now ready to answer the question we posed above regarding the
relationship between $P_{1,v}(r)/r^{n}$ and $\tilde{V}(t)$. (Compare to
Example~\ref{MCFdensity}.)

\begin{corollary}
\label{SolitonReducedVolume}Let $(\mathcal{M}^{n},g(t))$ be a compact
shrinking Ricci soliton that vanishes into a parabolic vertex $\mathcal{O}$ at
time $T$. Then for all $t<T$ and $r>0$, one has%
\[
\Theta_{\mathcal{O}}^{\mathrm{RF}}(t):=\tilde{V}_{\mathcal{O}}(t)=\frac
{P_{1,v}(r)}{r^{n}},
\]
where $P_{1,v}(r)=\int_{E_{r}}[|\nabla\ell|^{2}+R(n\log\frac{r}{\sqrt{4\pi
\tau}}-\ell)]\,d\mu\,dt$ is computed with $\ell=\ell_{\mathcal{O}}$.
\end{corollary}

\begin{proof}
It will be easiest to regard everything as a function of $\tau(t):=T-t>0$.
Because $(\mathcal{M}^{n},g(\tau))$ is a compact shrinking soliton, there
exist a time-independent metric $\bar{g}$ and function $\bar{f}$ on
$\mathcal{M}^{n}$ such that $\operatorname*{Rc}(\bar{g})+\bar{\nabla}%
\bar{\nabla}\bar{f}-\frac{1}{2}\bar{g}=0$. The solution of Ricci flow is then
$g(\tau)=\tau\xi_{\tau}^{\ast}(\bar{g})$, where $\{\xi_{\tau}\}_{\tau>0}$ is a
one-parameter family of diffeomorphisms such that $\xi_{1}=\operatorname*{id}$
and $\frac{\partial}{\partial\tau}\xi_{\tau}(x)=-\tau^{-1}\operatorname*{grad}%
\!_{\bar{g}}\bar{f}(x)$. The soliton potential function satisfies
$f(x,\tau)=\xi_{\tau}^{\ast}\bar{f}(x)$ and $f_{\tau}=-\left\vert \nabla
f\right\vert ^{2}$. (Notice that (\ref{RSE}) implies that system (\ref{RSS}) holds.)

Let $\Psi=(4\pi\tau)^{-n/2}e^{-\ell(x,\tau)}$, where $\ell$ is the reduced
distance from the parabolic vertex $\mathcal{O}$. By
Lemma~\ref{CLN-Observation}, $\ell=f+C$. So Assumptions~\ref{LL} and
\ref{Omega} are clearly satisfied. Because%
\[
\int_{\mathcal{M}^{n}\times\{\tau\}}|\psi|\,d\mu=O\mathcal{[}\tau^{n/2}%
\log(\tau^{-n/2})]\quad\text{and}\quad\int_{\mathcal{M}^{n}\times\{\tau
\}}|\nabla\psi|^{2}\,d\mu=O(\tau^{n/2-1})
\]
as $\tau\searrow0$, Assumption~\ref{r-bar} is satisfied as well. Because
$\frac{\partial}{\partial\tau}\psi=|\nabla f|^{2}-\frac{n}{2\tau}$,
Lemma~\ref{AlternateP} implies that%
\[
P_{1,v}(r)=\int_{E_{r}}(\frac{n}{2\tau}+\ell_{\tau}+|\nabla\ell|^{2}%
)\,d\mu\,dt\circeq\int_{E_{r}}\frac{n}{2\tau}\,d\mu\,dt.
\]
(Compare Remark~\ref{RD-alt}.) Computing $\tilde{V}(\tau)=\tilde
{V}_{\mathcal{O}}(\tau)$, one finds that%
\begin{align*}
\tilde{V}(1)  &  =\int_{\mathcal{M}^{n}}(4\pi)^{-n/2}e^{-\ell(x,1)}%
\,d\mu(g(1))\\
&  =\int_{0}^{\infty}\operatorname*{Vol}\!_{g(1)}\{x:\left(  4\pi\right)
^{-n/2}e^{-\ell\left(  x,1\right)  }\geq z\}\,dz\\
&  \circeq\int_{0}^{\infty}\frac{n}{2\tau}\operatorname*{Vol}\!_{g(\tau
)}\left[  \xi_{\tau}^{-1}\{x:(4\pi\tau)^{-n/2}e^{-\ell(x,1)}\geq1\}\right]
\,d\tau\qquad(z=\tau^{n/2})\\
&  =\int_{0}^{\infty}\frac{n}{2\tau}\operatorname*{Vol}\!_{g(\tau)}%
\{y:\ell(y,\tau)<n\log\frac{1}{\sqrt{4\pi\tau}}\}\,d\tau\\
&  =\int_{E_{1}}\frac{n}{2\tau}\,d\mu\,dt\\
&  \circeq P_{1,v}(1).
\end{align*}
But on a shrinking gradient soliton, $P_{1,v}(r)/r^{n}$ is independent of
$r>0$, while $\tilde{V}(\tau)$ is independent of $\tau>0$. Since they agree at
$r=1$ and $\tau=1$, they agree everywhere.
\end{proof}

Since the reduced distance and reduced volume are invariant under parabolic
rescaling, similar considerations apply to solutions whose rescaled limits are
shrinking gradient solitons.

\subsection{Localizing forward reduced volume}

In \cite{FIN05}, Feldman, Ilmanen, and the third author introduce a
\emph{forward reduced distance}%
\[
\ell_{+}(x,t):=\inf_{\gamma}\frac{1}{2\sqrt{t}}\int_{0}^{t}\sqrt{s}\left(
|\frac{d\gamma}{ds}|^{2}+R\right)  \,ds.
\]
Here the infimum is taken over smooth paths $\gamma$ from an origin $(\bar
{x},0)$ to $(x,t)$. Define%
\[
u(x,t)=(4\pi t)^{-n/2}e^{-\ell_{+}(x,t)}%
\]
and $\psi=\log u$. In \cite{Ni-LYH}, it is proved that $(\frac{\partial
}{\partial t}-\Delta-R)u\leq0$ holds in the distributional sense if
$(\mathcal{M}^{n},g(t))$ is a complete solution of Ricci flow with bounded
nonnegative curvature operator for $0\leq t\leq T$. Following the same
arguments as in the proof of Corollary~\ref{ReducedDistanceCorollary} then
leads to the following result for%
\[
P_{\varphi,u}(r)=\int_{E_{r}}[|\nabla\ell_{+}|^{2}-R(n\log\frac{r}{\sqrt{4\pi
t}}-\ell_{+})]\varphi\,d\mu\,dt.
\]

\begin{corollary}
\label{FRD}Let $(\mathcal{M}^{n},g(t))$ be a complete solution of Ricci flow
with bounded nonnegative curvature operator for $0\leq t\leq T$. Let $\varphi$
be any smooth nonnegative function. Then whenever $0<r_{0}<r_{1}<\sqrt{4\pi
T}$, one has%
\begin{equation}
\frac{P_{\varphi,u}(r_{1})}{r_{1}^{n}}-\frac{P_{\varphi,u}(r_{0})}{r_{0}^{n}%
}\leq\int_{r_{0}}^{r_{1}}\frac{n}{r^{n+1}}\int_{E_{r}}(\psi+n\log
r)(\frac{\partial\varphi}{\partial t}+\Delta\varphi)\,d\mu\,dt\,dr.
\label{ForwardVolumeComparison}%
\end{equation}

\end{corollary}

In direct analogy with Corollary~\ref{VolumeEquality}, one also has the following.

\begin{corollary}
Let $(\mathcal{M}^{n},g(t))$ be a complete solution of Ricci flow with bounded
nonnegative curvature operator for $0\leq t\leq T$. If equality holds in
(\ref{ForwardVolumeComparison}) with $\varphi\equiv1$, then $(E_{r},g(t))$ is
isometric to an expanding gradient soliton for all $r<\sqrt{4\pi T}$.
\end{corollary}

\section{Average energy for Ricci flow\label{AE}}

Again assume $(\mathcal{M}^{n},g(t))$ is a smooth complete solution of Ricci
flow for $t\in\lbrack0,\bar{t}]$. Let $\Psi$ denote a fundamental solution to
the conjugate heat equation
\begin{equation}
(\frac{\partial}{\partial t}+\Delta-R)\Psi=0 \label{CHE}%
\end{equation}
centered at $(\bar{x},\bar{t})$. The traditional notation in this case is
$\Psi=e^{-f}$, i.e.~$f:=-\psi$.

Perelman \cite{Perelman1} has discovered that the average energy%
\[
\mathcal{F}(t)=\int_{\mathcal{M}^{n}\times\{t\}}(\Delta f+R)e^{-f}\,d\mu
=\int_{\mathcal{M}^{n}\times\{t\}}(|\nabla f|^{2}+R)e^{-f}\,d\mu
\]
is a monotonically (weakly) increasing function of $t$. Our result in this
situation gives a quantity which is not just monotonic, but constant in its parameter.

\begin{corollary}
\label{Entropy}Suppose that $(\mathcal{M}^{n},g(t))$ is a smooth, compact
solution of Ricci flow for $t\in\lbrack0,\bar{t}]$, with $\bar{t}<\infty$.
Suppose further that $\Psi:\mathcal{M}^{n}\times\lbrack0,\bar{t}%
]\rightarrow(0,\infty)$ is a fundamental solution of (\ref{CHE}) with
singularity at $(\bar{x},\bar{t})$. Define $f:=-\log\Psi$.

Then for all $\bar{f}\in\mathbb{R}$ below some threshold value, we have
\[
\int_{\{f<\bar{f}\}}(\Delta f+R)e^{-\bar{f}}\,d\mu\,dt=1,
\]
where
\[
\{f<\bar{f}\}:=\{(x,t)\in\mathcal{M}^{n}\times\lbrack a,b):f(x,t)<\bar{f}\}.
\]

\end{corollary}

\begin{proof}
The arguments in Section~\ref{MVT} (below) verify that the hypotheses of
Lemma~\ref{AlternateP} are satisfied. Since
\[
\frac{\partial\psi}{\partial t}+|\nabla\psi|^{2}=-\Delta\psi+R,
\]
one then has%
\[
P_{1,\Psi}(r)=\int_{E_{r}}(-\Delta\psi+R)\,d\mu\,dt=\int_{E_{r}}(\Delta
f+R)\,d\mu\,dt.
\]
At this point, we change variables from $r$ to $\bar{f}:=n\log r$. We then
get
\[
\frac{P_{1,\Psi}(r)}{r^{n}}=\int_{\{f<\bar{f}\}}(\Delta f+R)e^{-\bar{f}}%
\,d\mu\,dt,
\]
whence the conclusion follows from Corollary~\ref{Moving} in Section~\ref{MVT}.
\end{proof}

\section{Mean-value theorems for heat kernels\label{MVT}}

In this section, we apply Theorem~\ref{MainTheorem} to heat kernels of
evolving Riemannian manifolds, especially those evolving by Ricci flow, with
stationary (i.e.~time-independent) manifolds appearing as an interesting
special case.

\medskip

Let $(\mathcal{M}^{n},g(t))$ be a smooth family of Riemannian manifolds
evolving by (\ref{GeneralMetricEvolution}) for $t\in\lbrack0,\bar{t}]$. We
will again abuse notation by regarding certain evolving quantities, where
convenient, as functions of $x\in\mathcal{M}^{n}$ and $\tau(t):=\bar{t}-t$.

A smooth function $\Psi:(\mathcal{M}^{n}\times\lbrack0,\bar{t}])\backslash
(\bar{x},0)\rightarrow\mathbb{R}_{+}$ is called a \emph{fundamental solution}
of the conjugate heat equation
\begin{equation}
(\frac{\partial}{\partial\tau}-\Delta-\operatorname*{tr}\!_{g}h)\Psi=0
\label{ConjugateHeat}%
\end{equation}
with singularity at $(\bar{x},0)$ if $\Psi$ satisfies (\ref{ConjugateHeat})
at$\ $all $(x,\tau)\in\mathcal{M}^{n}\times(0,\bar{t}]$, with $\lim
_{\tau\searrow0}\Psi(\cdot,\tau)=\delta_{\bar{x}}$ in the sense of
distributions. We call a minimal fundamental solution of (\ref{ConjugateHeat})
a \emph{heat kernel.}

For any smooth family $(\mathcal{M}^{n},g(t))$ of complete Riemannian
manifolds, it is well known that a heat kernel $\Psi$ always exists and is
unique. Moreover, $\Psi$ is bounded outside any compact space-time set
containing $(\bar{x},0)$ in its interior.\footnote{There are several standard
constructions, all of which utilize local properties that the manifold
inherits from $\mathbb{R}^{n}$. See the fine survey \cite{Gr99} and references
therein.} If $\Psi$ is the conjugate heat kernel for $(\mathcal{M}^{n},g(t))$,
then (\ref{DefineP}) takes the form
\[
P_{\varphi,\Psi}(r)=\int_{E_{r}}[|\nabla\log\Psi|^{2}-\log(r^{n}%
\Psi)(\operatorname*{tr}\!_{g}h)]\varphi\,d\mu\,dt.
\]
It is clear that Assumptions~\ref{LL} and ~\ref{Omega} are always satisfied.
In particular, $E_{\bar{r}}$ is compact for $\bar{r}>0$ sufficiently small. We
shall prove that Assumption~\ref{r-bar} is also valid for such $\bar{r}$. For
this, we need a purely local observation about $\Psi$ near $(\bar{x},0)$.

\begin{lemma}
\label{Friedman}For $t\in\lbrack0,\bar{t}]$, let $(\mathcal{N}^{n},g(t))$ be a
smooth family of (possibly incomplete) Riemannian manifolds. Suppose that
$\Psi$ is any fundamental solution of (\ref{ConjugateHeat}) with singularity
at $(\bar{x},0)$. For any $\varepsilon>0$, there exist a precompact
neighborhood $\Xi$ of $\bar{x}$, a time $\bar{\tau}\in(0,\bar{t}]$, and a
smooth function $\Phi:\Xi\times\lbrack0,\bar{\tau}]\rightarrow\mathbb{R}_{+}$
with $\Phi(\bar{x},0)=1$ such that for all $(x,\tau)\in(\Xi\times\lbrack
0,\bar{\tau}])\backslash(\bar{x},0)$, one has%
\begin{equation}
\left\vert \Psi(x,\tau)-\Phi(x,\tau)\cdot\frac{1}{(4\pi\tau)^{n/2}}\exp\left(
-\frac{d_{g(\tau)}^{2}(\bar{x},x)}{4\tau}\right)  \right\vert \leq\varepsilon.
\label{FriedmanEstimate}%
\end{equation}

\end{lemma}

\begin{proof}
One begins with Garofalo and Lanconelli's asymptotics \cite[Theorem~2.1]{GL89}
for a fundamental solution with respect to a Riemannian metric on
$\mathbb{R}^{n}$ which is Euclidean outside of an arbitrarily large compact
neighborhood of the origin. The first step is a straightforward adaptation of
their proof to the case $h\neq0$. The second step is to glue a large ball
centered at $\bar{x}\in\mathcal{N}^{n}$ into Euclidean space, obtaining a
manifold $(\mathbb{R}^{n},\tilde{g}(t))$ which is identical to $(\mathcal{N}%
^{n},g(t))$ on a large neighborhood of $\bar{x}$ and to which the refined
asymptotics apply. The difference of the fundamental solutions $\Psi$ and
$\tilde{\Psi}$ for $(\mathcal{N}^{n},g(t))$ and $(\mathbb{R}^{n},\tilde
{g}(t))$, respectively, starts at zero as a distribution. By the comparison
principle, it stays uniformly small for a short time.
\end{proof}

We now consider Assumption~\ref{r-bar}. Let $\bar{r}>0$ be given. Apply
Lemma~\ref{Friedman} with $\varepsilon=\bar{r}^{-n}/2$. By shrinking $\Xi$ and
$\bar{\tau}$ if necessary, we may assume without loss of generality that
$1/2\leq\Phi\leq2$ in $\Xi\times\lbrack0,\bar{\tau}]$. Because $\Psi
(\cdot,\tau)\rightarrow\delta_{\bar{x}}$ as $\tau\searrow0$, we may also
assume $\bar{\tau}>0$ is small enough that $E_{\bar{r}}(\tau)\subseteq\Xi$ for
all $\tau\in(0,\bar{\tau}]$, where $E_{\bar{r}}(\tau):=E_{\bar{r}}%
\cap(\mathcal{M}^{n}\times\{\tau\})$. Then in $\bigcup\nolimits_{\tau
\in(0,\bar{\tau}]}E_{\bar{r}}(\tau)$, one has%
\begin{equation}
\frac{1}{(4\pi\tau)^{n/2}}\exp\left(  -\frac{d_{g(\tau)}^{2}(\bar{x},x)}%
{4\tau}\right)  \geq\frac{\Psi(x,\tau)-\varepsilon}{\Phi(x,\tau)}\geq\frac
{1}{4\bar{r}^{n}}, \label{Gaussian}%
\end{equation}
which implies that $d_{g(\tau)}^{2}(\bar{x},\cdot)\leq4\tau\lbrack\frac{n}%
{2}\log\frac{1}{\tau}+\log4-\frac{n}{2}\log(4\pi)+\log\bar{r}^{n}]$ there.
Reduce $\bar{\tau}>0$ if necessary so that $\bar{\tau}\leq4^{(n-2)/n}\pi
\bar{r}^{-2}$. Then one has%
\begin{equation}
d_{g(\tau)}^{2}(\bar{x},\cdot)\leq4n\tau\log\frac{1}{\tau} \label{diameter}%
\end{equation}
in $E_{\bar{r}}(\tau)$ for all $\tau\in(0,\bar{\tau}]$. Since $\Psi>\bar
{r}^{-n}=2\varepsilon$ in $E_{\bar{r}}$, one also has%
\begin{equation}
\frac{\Psi}{2}\leq\Psi-\varepsilon\leq\frac{\Phi(x,\tau)}{(4\pi\tau)^{n/2}%
}\leq\frac{2}{(4\pi\tau)^{n/2}}. \label{Psi-estimate}%
\end{equation}
If necessary, reduce $\bar{\tau}>0$ further so $\bar{\tau}\leq\bar{r}^{-1}$
and $\bar{\tau}\leq4^{(n-2)/n}\pi$. Then $\psi:=\log\Psi$ satisfies%
\[
|\psi|\leq n\log\frac{1}{\tau}%
\]
in $E_{\bar{r}}(\tau)$ for all $\tau\in(0,\bar{\tau}]$. By (\ref{diameter}),
this proves that $\lim_{\tau\searrow0}\int_{E_{\bar{r}}(\tau)}|\psi|\,d\mu=0$.

If $\bar{\tau}\leq\bar{r}^{2}$, then $\Psi\leq\bar{r}^{-n}\leq\bar{\tau
}^{-n/2}$ outside $E_{\bar{r}}$. So by (\ref{Psi-estimate}), there exists
$c=c(n)$ such that $\Psi(\cdot,\tau)\leq e^{c}\tau^{-n/2}$ for all $\tau
\in(0,\bar{\tau}]$. By (\ref{diameter}), $E_{\bar{r}}(\tau)\subseteq
B_{g(\tau)}(\bar{x},\rho)$ for $\rho:=\sqrt{5n/e}$. Since $\Psi>\bar{r}^{-n}$
in $E_{\bar{r}}$, Theorem~\ref{LocalGradientEstimate} yields $C$ independent
of $x$ and $\tau$ in $\bigcup\nolimits_{0<\tau\leq\bar{\tau}}B_{g(\tau
)}\left(  \bar{x},2\rho\right)  $ such that for any $\tau\in(0,\bar{\tau}]$,
one has
\begin{equation}
|\nabla\psi|^{2}\leq(\frac{1}{\tau}+C)(1+c+n\log\bar{r}-\frac{n}{2}\log
\frac{\tau}{2})^{2} \label{gradient}%
\end{equation}
in $E_{\bar{r}}\cap(\mathcal{M}^{n}\times\lbrack\tau/2,\tau])$. If $\bar{r}>0$
is small enough that $E_{\bar{r}}$ is compact, this estimate and
(\ref{diameter}) prove that $\int_{E_{\bar{r}}}|\nabla\psi|^{2}\,d\mu
\,dt<\infty$, which establishes Assumption~\ref{r-bar}.

\begin{remark}
Assumption~\ref{r-bar} is valid for \emph{all} $\bar{r}>0$ in any manifold
$(\mathcal{M}^{n},g(t))_{t<\bar{t}}$ for which the kernel $\Psi$ vanishes at
infinity in space-time, i.e.~if for every $\varepsilon>0$, there exists a
compact set $K\subset\mathcal{M}^{n}\times(-\infty,\bar{t}]$ such that
$\Psi\leq\varepsilon$ outside $K$.
\end{remark}

Our main result in this section is the following consequence of
Theorem~\ref{MainTheorem}. The reader is invited to compare it with
Corollary~\ref{ReducedDistanceCorollary} (above) for Perelman's reduced volume
density. Recall that $\psi_{(r)}:=\log(r^{n}\Psi)$.

\begin{corollary}
\label{Moving}Suppose that $(\mathcal{M}^{n},g(t))$ is a smooth family of
complete Riemannian manifolds evolving by (\ref{GeneralMetricEvolution}) for
$t\in\lbrack0,\bar{t}]$. Let $\Psi:(\mathcal{M}^{n}\times\lbrack0,\bar
{t}])\backslash(\bar{x},0)\rightarrow\mathbb{R}_{+}$ be the kernel of the
conjugate heat equation (\ref{ConjugateHeat}) with singularity at
$(x,\tau)=(\bar{x},0)$. Let $\varphi$ be any smooth function of $(x,t)$. Then
there is $\bar{r}>0$ such that if $0<r_{0}<r_{1}<\bar{r}$, then%
\[
\frac{P_{\varphi,\Psi}(r_{1})}{r_{1}^{n}}-\frac{P_{\varphi,\Psi}(r_{0})}%
{r_{0}^{n}}=-\int_{r_{0}}^{r_{1}}\frac{n}{r^{n+1}}\int_{E_{r}}\psi_{(r)}%
(\frac{\partial\varphi}{\partial t}-\Delta\varphi)]\,d\mu\,dt\,dr.
\]
Furthermore, one has%
\[
\varphi(\bar{x},\bar{t})=\lim_{r\searrow0}\frac{P_{\varphi,\Psi}(r)}{r^{n}},
\]
and thus%
\[
\varphi(\bar{x},\bar{t})=\frac{P_{\varphi,\Psi}(r_{1})}{r_{1}^{n}}+\int
_{0}^{r_{1}}\frac{n}{r^{n+1}}\int_{E_{r}}\psi_{(r)}(\frac{\partial\varphi
}{\partial t}-\Delta\varphi)]\,d\mu\,dt\,dr.
\]

\end{corollary}

\begin{proof}
Now that we have verified Assumptions~\ref{LL}--\ref{r-bar}, everything
follows directly from Theorem~\ref{MainTheorem} except for the representation
formula $\varphi(\bar{x},\bar{t})=\lim_{r\searrow0}[P_{\varphi,\Psi}%
(r)/r^{n}]$, which we will prove by a blow-up argument. Without loss of
generality, we may assume that $\varphi(\bar{x},\bar{t})=1$. Here is the
set-up. Identify $\mathbb{R}^{n}$ with $T_{\bar{x}}\mathcal{M}^{n}$, and let
$y\in\mathcal{M}^{n}$ denote the image of $\vec{y}\in\mathbb{R}^{n}$ under the
exponential map $\exp_{\bar{x}}(\cdot)$ for $g$ at $\tau=0$. For $r>0$, define
$\varphi^{r}(\vec{y},\tau):=\varphi(ry,r^{2}\tau)$, $\Psi^{r}(\vec{y}%
,\tau):=r^{n}\Psi(ry,r^{2}\tau)$, and $\Psi^{0}(\vec{y},\tau):=(4\pi
\tau)^{-n/2}e^{-|\vec{y}|^{2}/4\tau}$. Let $d\mu^{r}(\cdot,\tau)$ denote the
pullback of $r^{-n}d\mu(\cdot,r^{2}\tau)$ under the map $\vec{y}\mapsto
\exp_{\bar{x}}(r\vec{y})$. For $\delta\geq0$, consider the `truncations'
defined by
\begin{align*}
E_{r}^{\delta}  &  :=E_{r}\cap(\mathcal{M}^{n}\times(\delta r^{2},\bar
{t}\,]),\\
\hat{E}_{r}^{\delta}  &  :=\{(\vec{y},\tau):\tau>\delta\text{ and }\Psi
^{r}(\vec{y},\tau)>1\},\\
\hat{E}_{0}^{\delta}  &  :=\{(\vec{y},\tau):\tau>\delta\text{ and }\Psi
^{0}(\vec{y},\tau)>1\},\\
P_{r}^{\delta}  &  :=\int_{E_{r}^{\delta}}[|\nabla\log\Psi|^{2}%
-(\operatorname*{tr}\!_{g}h)\log(r^{n}\Psi)]\varphi\,d\mu\,d\tau\\
P_{0}^{\delta}  &  :=\int_{\hat{E}_{0}^{\delta}}|\nabla\log\Psi^{0}%
|^{2}\,d\vec{y}\,d\tau.
\end{align*}
The proof consists of two claims, which together imply the result.

The first claim is that if $0<\delta\ll1$, then $\lim_{r\searrow0}%
[P_{r}^{\delta}/r^{n}]=P_{0}^{\delta}$. Pulling back, one computes
$P_{r}^{\delta}=r^{n}\int_{\hat{E}_{r}^{\delta}}[|\nabla\log\Psi^{r}%
|^{2}-r^{2}(\operatorname*{tr}\!_{g}h)\log\Psi^{r}]\varphi^{r}\,d\mu
^{r}\,d\tau$. By Lemma~\ref{Friedman}, $\Psi^{r}\rightarrow\Psi^{0}$ as
$r\searrow0$ uniformly on any $\Omega\subset\subset\mathbb{R}^{n}\times
\lbrack\delta,\bar{t}]$. By parabolic regularity, $\chi(\hat{E}_{r}^{\delta
})\rightarrow\chi(\hat{E}_{0}^{\delta})$ in $L^{1}(\mathbb{R}^{n})$ as
$r\searrow0$. Since $d\mu^{r}\rightarrow d\vec{y}$ and $\varphi(\bar{x}%
,\bar{t})=1$, the claim follows.

The second claim is that for any $\eta>0$, there exists some $\delta
\in(0,1/100)$ such that $0\leq\lbrack P_{\varphi,\Psi}(r)-P_{r}^{\delta
}]/r^{n}<\eta$ for all small $r>0$. By Lemma~\ref{Friedman}, if $r\leq1$ is so
small that $1/2\leq\Phi\leq2$ in $E_{r}\cap(\Xi\times\lbrack0,\bar{\tau}])$,
then $(4\pi\tau)^{-n/2}\exp\left(  -d_{g(\tau)}^{2}(\bar{x},x)/4\tau\right)
\geq\frac{1}{4r^{n}}$ there. (Compare (\ref{Gaussian}).) Furthermore,
$d_{g(\tau)}^{2}(\bar{x},\cdot)\leq4\tau(\frac{n}{2}\log\frac{r^{2}}{\tau
}+\log4)\leq4\tau n\log\frac{r^{2}}{\tau}$ in $E_{r}\backslash E_{r}^{\delta}%
$, since $\frac{r^{2}}{\tau}\geq4$. Because $\tau\leq1$ in $E_{r}$ for all
small $r>0$, Theorem~\ref{LocalGradientEstimate} gives $C$ such that
$|\nabla\log\Psi|^{2}\leq\frac{C}{\tau}(\log\frac{r^{2}}{\tau})^{2}$ in
$E_{r}\backslash E_{r}^{\delta}$. (Here we used $r^{-n}\leq\Psi\leq e^{c}%
\tau^{-n/2}$; compare (\ref{gradient}).) Therefore,%
\[
\int_{E_{r}\backslash E_{r}^{\delta}}|\nabla\log\Psi|^{2}\,d\mu\leq C^{\prime
}\int_{0}^{\delta r^{2}}\tau^{\frac{n-2}{2}}\left(  \log\frac{r^{2}}{\tau
}\right)  ^{\frac{n+4}{2}}d\tau\leq C^{\prime\prime}r^{n}\delta^{n/2}\left(
\log\frac{1}{\delta}\right)  ^{\frac{n+4}{2}}.
\]
The second claim, hence the theorem, follows readily.
\end{proof}

\begin{remark}
In the special case that $L(\cdot,t)$ is a divergence-form, uniformly elliptic
operator on Euclidean space $\mathbb{R}^{n}$ and $\Psi$ is the kernel of its
adjoint $L^{\ast}$, the results of Corollary~\ref{Moving} appear in
\cite[Theorems~1 and 2]{FG87} for $\varphi$ solving $(\frac{\partial}{\partial
t}-L)\varphi=0$, and in \cite[Theorem~1.5]{GL88} for arbitrary smooth
$\varphi$.
\end{remark}

\medskip

We conclude this section with two results for the special case of the
conjugate heat kernel $\Psi$ of a \emph{fixed} Riemannian manifold
$(\mathcal{M}^{n},g)$.

\bigskip

Our first observation is that one can adapt the argument of \cite{FG87} to
obtain a mean-value representation theorem in terms of an integral on `heat
spheres'. This approach is naturally related to the interpretation of equation
(\ref{STG2-1}) as a space-time Green's formula. To give the argument, we
introduce some additional notation. Consider the space-time manifold
$\widetilde{\mathcal{M}}^{n+1}=\mathcal{M}^{n}\times\mathbb{R}$ equipped with
the metric $\tilde{g}(x,t)=g(x)+dt^{2}$, where $t$ is the global $\mathbb{R}%
$-coordinate. Applying Green's formula to a bounded space-time domain $D$ in
$\widetilde{\mathcal{M}}^{n+1}$ with the vector field $\varphi\Psi
\frac{\partial}{\partial t}-\Psi\nabla\varphi+\varphi\nabla\Psi$, we get%
\begin{align}
\int_{D}(\frac{\partial\varphi}{\partial t}-\Delta\varphi)\Psi\,d\mu\,dt  &
=\int_{D}[(\frac{\partial\varphi}{\partial t}-\Delta\varphi)\Psi+\varphi
(\frac{\partial\Psi}{\partial t}+\Delta\Psi)]\,d\mu\,dt\label{STDivergence}\\
&  =\int_{D}\operatorname*{div}\!_{\tilde{g}}(\varphi\Psi\frac{\partial
}{\partial t}-\Psi\nabla\varphi+\varphi\nabla\Psi)\,d\mu\,dt\nonumber\\
&  =\int_{\partial D}\left\langle \varphi\Psi\frac{\partial}{\partial t}%
-\Psi\nabla\varphi+\varphi\nabla\Psi,\,\tilde{\nu}\right\rangle _{\tilde{g}%
}\,d\tilde{A},\nonumber
\end{align}
where $\tilde{\nu}$ is the unit outward normal and $d\tilde{A}$ the area
element of $\partial D$, both taken with respect to $\tilde{g}$. For $s\geq0$,
we follow \cite{FG87} in defining%
\[
D_{r}^{s}=\{(x,\tau)\in E_{r}:\tau>s\}
\]
and two portions of its space-time boundary,%
\[
P_{1}^{s}=\{(x,\tau):\Psi=r^{-n},\;\tau>s\}\qquad\text{and}\qquad P_{2}%
^{s}=\{(x,\tau)\in\bar{D}_{r}^{s}:\tau=s\}.
\]
Applying (\ref{STDivergence}) to $D_{r}^{s}$ yields%
\begin{align*}
0  &  =\int_{D_{r}^{s}}(\frac{\partial\varphi}{\partial t}-\Delta\varphi
)\Psi\,d\mu\,dt\\
&  =\int_{P_{2}^{s}}\varphi\Psi\,d\mu+\int_{P_{1}^{s}}\left\langle \varphi
\Psi\frac{\partial}{\partial t}-\Psi\nabla\varphi+\varphi\nabla\Psi
,\,\tilde{\nu}\right\rangle _{\tilde{g}}\,d\tilde{A}\\
&  =\int_{P_{2}^{s}}\varphi\Psi\,d\mu+\frac{1}{r^{n}}\int_{P_{1}^{s}%
}\left\langle \varphi\frac{\partial}{\partial t}-\nabla\varphi,\,\tilde{\nu
}\right\rangle _{\tilde{g}}\,d\tilde{A}+\int_{P_{1}^{s}}\varphi\left\langle
\nabla\Psi,\,\tilde{\nu}\right\rangle _{\tilde{g}}\,d\tilde{A}.
\end{align*}
Letting $s\searrow0$, we obtain%
\begin{align*}
\varphi(\bar{x},0)  &  =\lim_{s\searrow0}\int_{P_{2}^{s}}\varphi\Psi\,d\mu\\
&  =-\frac{1}{r^{n}}\int_{P_{1}^{0}}\left\langle \varphi\frac{\partial
}{\partial t}-\nabla\varphi,\,\tilde{\nu}\right\rangle _{\tilde{g}}%
\,d\tilde{A}-\int_{P_{1}^{0}}\varphi\left\langle \nabla\Psi,\,\tilde{\nu
}\right\rangle _{\tilde{g}}\,d\tilde{A}\\
&  =-\frac{1}{r^{n}}\int_{D_{r}^{0}}(\frac{\partial}{\partial t}%
-\Delta)\varphi\,d\mu\,dt+\int_{P_{1}^{0}}\varphi\frac{|\nabla\Psi|^{2}}%
{\sqrt{|\Psi_{t}|^{2}+|\nabla\Psi|^{2}}}\,d\tilde{A}.
\end{align*}
Summing together and noticing that $P_{1}^{0}=\partial E_{r}$, we get the
following mean-value theorem, which is naturally related to
Corollary~\ref{Moving} by the coarea formula.

\begin{theorem}
\label{LeiNi2}Let $(\mathcal{M}^{n},g)$ be a complete fixed manifold. Let
$\Psi$ denote the conjugate heat kernel with singularity at $(x,\tau)=(\bar
{x},0)$. If a smooth function $\varphi$ of $(x,t)$ solves the heat equation,
then%
\[
\varphi(\bar{x},\bar{t})=\int_{\partial E_{r}}\frac{|\nabla\Psi|^{2}}%
{\sqrt{|\Psi_{t}|^{2}+|\nabla\Psi|^{2}}}\varphi\,d\tilde{A}.
\]

\end{theorem}

\medskip

For the $\varepsilon$-regularity theorems for Ricci flow derived by the third
author \cite{Ni-MVT}, we need a mean-value inequality for nonnegative
supersolutions. For this purpose, assume that the Ricci curvature of
$(\mathcal{M}^{n},g)$ satisfies $\operatorname*{Rc}\geq(n-1)kg$ for some
$k\in\{-1,0,1\}$. Let $(\mathcal{M}_{k}^{n},\tilde{g})$ denote the simply
connected space form of constant sectional curvature $k$, and let $\Psi_{k}$
denote its conjugate heat kernel centered at $\tilde{x}\in\mathcal{M}_{k}^{n}%
$. Then there exists $\tilde{\Psi}_{k}:[0,\infty)\times(0,\infty
)\rightarrow(0,\infty)$ such that $\Psi_{k}(x,\tau)=\tilde{\Psi}_{k}%
(d_{k}(\tilde{x},x),\tau)$, where $d_{k}$ denotes the distance function of
$(\mathcal{M}_{k}^{n},\tilde{g})$.

Fix an origin $(\bar{x},\bar{\tau})\in\mathcal{M}^{n}\times\mathbb{R}$. Again
let $\tau:=\bar{t}-t$, and let $\tilde{\Psi}$ denote the transplant of
$\Psi_{k}$ to $(\mathcal{M}^{n},g)$, i.e.%
\begin{equation}
\tilde{\Psi}(x,\tau):=\tilde{\Psi}_{k}(d_{g}(\bar{x},x),\tau).
\label{Transplant}%
\end{equation}
As above, let $\tilde{\psi}_{(r)}=\log(r^{n}\tilde{\Psi})$ and $\tilde{E}%
_{r}=\{(x,\tau)\in\mathcal{M}^{n}\times\mathbb{R}:\tilde{\psi}_{(r)}%
(x,\tau)>0\}$. Define%
\begin{equation}
\tilde{I}^{(\bar{x},0)}(r):=\frac{1}{r^{n}}\int_{\tilde{E}_{r}}|\nabla
\log\tilde{\Psi}|^{2}\,d\mu\,dt.
\end{equation}
Then the following mean-value inequality follows from
Theorem~\ref{MainTheorem}.

\begin{corollary}
\label{LN-MVI}Let $(\mathcal{M}^{n},g)$ be a complete Riemannian manifold such
that $\operatorname*{Rc}\geq(n-1)kg$ for some $k\in\{-1,0,1\}$. Let
$\tilde{\Psi}$ be defined by (\ref{Transplant}), and let $\varphi\geq0$ be any
smooth supersolution of the heat equation, i.e.~$(\frac{\partial}{\partial
t}-\Delta)\varphi\geq0$. Then%
\[
\varphi(\bar{x},\bar{t})\geq\frac{1}{r^{n}}\int_{\tilde{E}_{r}}|\nabla
\log\tilde{\Psi}|^{2}\varphi\,d\mu\,dt.
\]
In particular, $\tilde{I}^{(\bar{x},0)}(r)\leq1$ holds for all $r>0$, and
$\frac{d}{dr}\tilde{I}^{(\bar{x},0)}(r)\leq0$ holds in the sense of distributions.

If equality holds for $\varphi\equiv1$, then the largest metric ball in
$\tilde{E}_{r}$ is isometric to the corresponding ball in the simply-connected
space form of constant sectional curvature $k$.
\end{corollary}

\begin{proof}
The inequalities follow from Theorem~\ref{MainTheorem} by the results of
Cheeger--Yau \cite{CY} that $(\frac{\partial}{\partial\tau}-\Delta)\tilde
{\Psi}(x,\tau)\leq0$ and $\tilde{\Psi}(x,\tau)\geq\Psi(x,\tau)$, where $\Psi$
is the conjugate heat kernel of $(\mathcal{M}^{n},g)$. The implication of
equality is a consequence of the rigidity derived from equality in the Bishop
volume comparison theorem. (See \cite{CE}.)
\end{proof}

\section{Appendix: simple estimates for reduced geometry\label{Appendix}}

For the convenience of the reader, we provide certain elementary estimates
involving reduced geometry in a form adapted to this paper. The reader should
note that most of the estimates solely for reduced distance are essentially
contained in Ye's notes \cite{Ye}, though not always in the form stated here.
(Also see \cite{RFV2}.)

\medskip

\noindent\textbf{Notation.} Assume that $(\mathcal{M}^{n},g(\tau))$ is a
smooth one-parameter family of complete (possibly noncompact) manifolds
satisfying $\frac{\partial}{\partial\tau}g=2\operatorname*{Rc}$ for $0\leq
\tau\leq\bar{\tau}$. Unless otherwise noted, all Riemannian quantities are
measured with respect to $g(\tau)$. All quantities in reduced geometry are
calculated with respect to a fixed origin $\mathcal{O}=(\bar{x},0)$. We denote
the metric distance from $x$ to $y$ with respect to $g(\tau)$ by $d_{\tau
}(x,y)$ and write $d_{\tau}(x)=d_{\tau}(\bar{x},x)$. We define $B_{\tau
}(x,r)=\{y\in\mathcal{M}^{n}:d_{\tau}(x,y)<r)$ and write $B_{\tau}(r)=B_{\tau
}(\bar{x},r)$. Perelman's \emph{space-time action} $\mathcal{L}$,
\emph{reduced distance} $\ell$, and \emph{reduced volume density} $v$ are
defined above in (\ref{ScriptL}), (\ref{ReducedDistance}), and
(\ref{ReducedVolumeDensity}), respectively. We will also use the
\emph{space-time distance} $L(x,\tau):=\inf\{\mathcal{L}(\gamma):\gamma
(0)=(\bar{x},0),~\gamma(\tau)=(x,\tau)\}$.

\subsection{Bounds for reduced distance}

Given $k\geq0$ and $K\geq0$, define%
\begin{equation}
\underline{\ell}(x,\tau)=e^{-2k\tau}\frac{d_{0}^{2}(x)}{4\tau}-\frac{nk}%
{3}\tau
\end{equation}
and%
\begin{equation}
\overline{\ell}(x,\tau)=e^{2K\tau}\frac{d_{0}^{2}(x)}{4\tau}+\frac{nK}{3}\tau.
\end{equation}
Our first observation directly follows Ye \cite{Ye}.

\begin{lemma}
\label{BoundsForReducedDistance}The reduced distance $\ell(x,\tau)$ has the
following properties.

\begin{enumerate}
\item If there is $k\geq0$ such that $\operatorname*{Rc}\geq-kg$ on
$\mathcal{M}^{n}\times\lbrack0,\bar{\tau}]$, then $\ell(x,\tau)\geq
\underline{\ell}(x,\tau)$.

\item If there is $K\geq0$ such that $\operatorname*{Rc}\leq Kg$ on
$\mathcal{M}^{n}\times\lbrack0,\bar{\tau}]$, then $\ell(x,\tau)\leq
\overline{\ell}(x,\tau)$.
\end{enumerate}
\end{lemma}

\begin{proof}
\ 

\begin{enumerate}
\item Observe that $g(\tau)\geq e^{-2k\tau}g(0)$. By (\ref{ScriptL}), the
$\mathcal{L}$-action of an arbitrary path $\gamma$ from $(\bar{x},0)$ to
$(x,\tau)$ is%
\begin{align*}
\mathcal{L}(\gamma)  &  =\int_{0}^{\sqrt{\tau}}\left(  \frac{1}{2}%
|\frac{d\gamma}{ds}|^{2}+2s^{2}R\right)  \,ds\\
&  \geq\frac{1}{2}e^{-2k\tau}\int_{0}^{\sqrt{\tau}}|\frac{d\gamma}{ds}%
|_{0}^{2}\,ds-2nk\int_{0}^{\sqrt{\tau}}s^{2}\,ds\\
&  \geq e^{-2k\tau}\frac{d_{0}^{2}(x)}{2\sqrt{\tau}}-\frac{2nk}{3}\tau^{3/2}.
\end{align*}
Since $\gamma$ was arbitrary, one has $\ell(x,\tau)=\frac{1}{2\sqrt{\tau}}%
\inf_{\gamma}\mathcal{L}(\gamma)\geq\underline{\ell}(x,\tau)$.

\item Observe that $g(\tau)\leq e^{2K\tau}g(0)$. Let $\beta$ be a path from
$(\bar{x},0)$ to $(x,\tau)$ that is minimal and of constant speed with respect
to $g(0)$. Then as above,
\[
\mathcal{L}(\beta)\leq e^{2K\tau}\frac{d_{0}^{2}(x)}{2\sqrt{\tau}}+\frac
{2nK}{3}\tau^{3/2}.
\]
Hence $\ell(x,\tau)\leq\frac{1}{2\sqrt{\tau}}\mathcal{L}(\beta)\leq
\overline{\ell}(x,\tau)$.
\end{enumerate}
\end{proof}

\begin{remark}
If $\operatorname*{Rc}\geq-kg$ on $\mathcal{M}^{n}\times\lbrack0,\bar{\tau}]$,
it follows from Part (1) of Lemma \ref{BoundsForReducedDistance} (by standard
arguments) that minimizing $\mathcal{L}$-geodesics exist and are smooth.
\end{remark}

\subsection{Bounds for reduced-volume heatballs}

Recall that the reduced-volume density is $v(x,\tau)=(4\pi\tau)^{-n/2}%
e^{-\ell(x,\tau)}$. For $r>0$, define the reduced-volume heatball%
\begin{align}
E_{r}  &  =\{(x,\tau)\in\mathcal{M}^{n}\times(0,\bar{\tau}]:v(x,\tau
)>r^{-n}\}\\
&  =\{(x,\tau)\in\mathcal{M}^{n}\times(0,\bar{\tau}]:\ell(x,\tau)<n\log
\frac{r}{\sqrt{4\pi\tau}}\}
\end{align}
and define $c(k,\bar{\tau})$ by
\begin{equation}
c=\frac{e^{4k\bar{\tau}/3}}{4\pi}. \label{DefineA}%
\end{equation}
Given $r>0$, $k\geq0$, $\tau>0$, define%
\begin{equation}
\rho(r,k,\tau)=e^{k\tau}\sqrt{(2n\tau\log\frac{r^{2}}{4\pi\tau}+\frac{4}%
{3}nk\tau^{2})_{+}}. \label{Define-rho}%
\end{equation}
Note that $\rho(r,0,\tau)$ agrees with $R_{r}(\tau)$ in \cite{EG82}. It is
easy to see that for each $r>0$ and $k\geq0$, one has $\rho(r,k,\tau)>0$ for
all sufficiently small $\tau>0$.

\begin{remark}
If $\operatorname*{Rc}\geq-kg$ on $\mathcal{M}^{n}\times\lbrack0,\bar{\tau}]$,
then Part (1) of Lemma \ref{BoundsForReducedDistance} implies that
$(x,\tau)\in E_{r}$ only if $x\in B_{0}(\rho(r,k,\tau))$.
\end{remark}

\begin{lemma}
\label{rhoEstimates}Assume that $0<r^{2}\leq\min\{\bar{\tau}/c,4\pi\}$. If
$cr^{2}\leq\tau\leq\bar{\tau}$, then $\rho(r,\tau)=0$.
\end{lemma}

\begin{proof}
When $\tau=cr^{2}$, one has%
\[
\frac{k}{3}\tau+\frac{1}{2}\log\frac{r^{2}}{4\pi\tau}\leq\frac{k\bar{\tau}}%
{3}+\frac{1}{2}\log\frac{1}{4\pi c}=-\frac{k\bar{\tau}}{3}\leq0,
\]
while for $cr^{2}\leq\tau\leq\bar{\tau}$, one has%
\begin{align*}
\frac{\partial}{\partial\tau}(\frac{k}{3}\tau^{2}+\frac{1}{2}\tau\log
\frac{r^{2}}{4\pi\tau})  &  =\frac{2k}{3}\tau+\frac{1}{2}\log\frac{r^{2}}%
{4\pi\tau}-\frac{1}{2}\\
&  \leq\frac{2k}{3}\bar{\tau}+\frac{1}{2}\log\frac{1}{4\pi c}-\frac{1}{2}%
\leq-\frac{1}{2}.
\end{align*}

\end{proof}

\begin{corollary}
\label{BoundReducedBalls}Assume that $\operatorname*{Rc}\geq-kg$ on
$\mathcal{M}^{n}\times\lbrack0,\bar{\tau}]$ for some $k\geq0$ and that
$0<r^{2}\leq\min\{\bar{\tau}/c,4\pi\}$. Then%
\[
E_{r}\subseteq%
{\displaystyle\bigcup\limits_{0<\tau<cr^{2}}}
B_{0}(\rho(r,k,\tau))\times\{\tau\}.
\]

\end{corollary}

\subsection{Gradient estimates for reduced distance\label{GED}}

Local gradient estimates for curvatures evolving by Ricci flow originated in
\cite[\S 7]{Shi89}. Recall the following version.

\begin{proposition}
[{Hamilton \cite[\S 13]{Ham95}}]\label{CurvatureGradientEstimate}Suppose
$g(\tau)$ solves backward Ricci flow for $\tau_{0}\leq\tau\leq\tau_{1}$ on an
open set $\mathcal{U}$ of $\mathcal{M}^{n}$ with $\bar{B}_{\tau_{1}%
}(x,2\lambda)\subset\mathcal{U}$. There exists $C_{n}$ depending only on $n$
such that if $|\operatorname*{Rm}|\leq M$ on $\mathcal{U}\times\lbrack\tau
_{0},\tau_{1}]$, then%
\[
|\nabla\operatorname*{Rm}|\leq C_{n}M\sqrt{\frac{1}{\lambda^{2}}+\frac{1}%
{\tau_{1}-\tau}+M}%
\]
on $B_{\tau_{1}}(x,\lambda)\times\lbrack\tau_{0},\tau_{1})$.
\end{proposition}

If there is a global bound on curvature, the situation is quite simple:

\begin{remark}
If $|\operatorname*{Rm}|\leq M$ on $\mathcal{M}^{n}\times\lbrack0,\bar{\tau}%
]$, then for every $\tau^{\ast}<\bar{\tau}$ there exists $A=A(n,M,\tau^{\ast
})$ such that $|\nabla R|\leq A$ on $\mathcal{M}^{n}\times\lbrack0,\tau^{\ast
}]$.
\end{remark}

More generally, the following `localization lemma' often provides adequate
local bounds.

\begin{lemma}
\label{Neighborhood}Assume $\operatorname*{Rc}\geq-kg$ on $\mathcal{M}%
^{n}\times\lbrack0,\bar{\tau}]$. Then for every $\lambda>0$ and $\tau^{\ast
}\in(0,\bar{\tau})$, there exists $\lambda^{\ast}$ such that the image of any
minimizing $\mathcal{L}$-geodesic from $(\bar{x},0)$ to any $(x,\tau)\in
B_{0}(\lambda)\times(0,\tau^{\ast}]$ is contained in $B_{0}(\lambda^{\ast})$.
Furthermore, there exist constants $C,C^{\prime}$ such that
$\operatorname*{Rc}<Cg$ and $|\nabla R|\leq C^{\prime}$ on $B_{0}%
(\lambda^{\ast})\times\lbrack0,\tau^{\ast}]$.
\end{lemma}

\begin{proof}
By smoothness, there exists $K$ such that $\operatorname*{Rc}\leq Kg$ on
$B_{0}(\lambda)\times\lbrack0,\tau^{\ast}]$. Applying Part (2) of Lemma
\ref{BoundsForReducedDistance} along radial geodesics from $\bar{x}$ shows
that
\[
\sup_{(x,\tau)\in B_{0}(\lambda)\times\lbrack0,\tau^{\ast}]}[\tau\ell
(x,\tau)]\leq e^{2K\tau^{\ast}}\frac{\lambda^{2}}{4}+\frac{nK}{3}(\tau^{\ast
})^{2}.
\]
Define%
\[
\lambda^{\ast}=2e^{k\tau^{\ast}}\sqrt{e^{2K\tau^{\ast}}\frac{\lambda^{2}}%
{4}+\frac{n}{3}(k+K)(\tau^{\ast})^{2}}.
\]
Let $(x,\tau)\in B_{0}(\lambda)\times(0,\tau^{\ast}]$ be arbitrary and let
$\gamma$ be any minimizing $\mathcal{L}$-geodesic from $(\bar{x},0)$ to
$(x,\tau)$. Then for every $\sigma\in\lbrack0,\tau]$, one obtains%
\[
d_{0}(\gamma(\sigma))\leq2e^{k\tau}\sqrt{\tau\ell(x,\tau)+\frac{nk}{3}\tau
^{2}}<\lambda^{\ast}%
\]
by following the proof of Part (1) of Lemma \ref{BoundsForReducedDistance}.
This proves that the image of $\gamma$ is contained in $B_{0}(\lambda^{\ast})$.

Now define $\tau^{\prime}=\tau^{\ast}+\frac{1}{2}(\bar{\tau}-\tau^{\ast})$ and
choose $\lambda^{\prime}$ large enough that $B_{0}(\lambda^{\ast})\subseteq
B_{\tau^{\prime}}(\lambda^{\prime})$. By smoothness, there exists $M$ such
that $|\operatorname*{Rm}|\leq M$ on $B_{\tau^{\prime}}(3\lambda^{\prime
})\times\lbrack0,\tau^{\prime}]$. So by Proposition
\ref{CurvatureGradientEstimate}, there exists $C^{\prime}$ such that $|\nabla
R|\leq C^{\prime}$ on $B_{\tau^{\prime}}(\lambda^{\prime})\times\lbrack
0,\tau^{\ast}]$. Clearly, $\operatorname*{Rc}\leq Cg$ on $B_{\tau^{\prime}%
}(\lambda^{\prime})\times\lbrack0,\tau^{\ast}]$ as well.
\end{proof}

\begin{lemma}
\label{GeodesicSpeed}Assume that there exists an open set $\mathcal{U}%
\subseteq\mathcal{M}^{n}$ and $0\leq\tau_{0}\leq\tau_{1}\leq\bar{\tau}$ such
that $|\nabla R|\leq A$ on $\mathcal{U}\times\lbrack\tau_{0},\tau_{1}]$. Let
$\gamma:[0,\tau_{1}]\rightarrow\mathcal{U}$ be an $\mathcal{L}$-geodesic and
let $\Gamma(\tau_{0})=\lim_{\tau\searrow\tau_{0}}(\sqrt{\tau}|\frac{d\gamma
}{d\tau}|)$, which is well defined for all $\tau_{0}\geq0$.

\begin{enumerate}
\item If $\operatorname*{Rc}\geq-kg$ on $\mathcal{U}\times\lbrack\tau_{0}%
,\tau_{1}]$, then for all $\tau\in\lbrack\tau_{0},\tau_{1}]$, one has
\begin{align*}
|\frac{d\gamma}{d\tau}|  &  \leq\frac{1}{2\sqrt{\tau}}\left[  (2\Gamma
(\tau_{0})+\frac{A}{k}\sqrt{\tau_{1}})e^{k(\tau-\tau_{0})}-\frac{A}{k}%
\sqrt{\tau_{1}}\right]  \qquad\;\;(k>0)\\
&  \leq\frac{1}{2\sqrt{\tau}}\left[  2\Gamma(\tau_{0})+A\sqrt{\tau_{1}}%
(\tau-\tau_{0})\right]  \qquad\qquad\qquad\qquad(k=0)
\end{align*}

\item If $\operatorname*{Rc}\leq Kg$ on $\mathcal{U}\times\lbrack\tau_{0}%
,\tau_{1}]$, then for all $\tau\in\lbrack\tau_{0},\tau_{1}]$, one has%
\begin{align*}
|\frac{d\gamma}{d\tau}|  &  \geq\frac{1}{2\sqrt{\tau}}\left[  (2\Gamma
(\tau_{0})+\frac{A}{K}\sqrt{\tau_{1}})e^{K(\tau_{0}-\tau)}-\frac{A}{K}%
\sqrt{\tau_{1}}\right]  \qquad\,(K>0)\\
&  \geq\frac{1}{2\sqrt{\tau}}\left[  2\Gamma(\tau_{0})+A\sqrt{\tau_{1}}%
(\tau_{0}-\tau)\right]  \qquad\qquad\qquad\qquad(K=0).
\end{align*}

\end{enumerate}
\end{lemma}

\begin{proof}
It will be more convenient to regard $\gamma$ as a function of $s=\sqrt{\tau}%
$. Let $\dot{\gamma}=\frac{d\gamma}{d\tau}$ and $\gamma^{\prime}=\frac
{d\gamma}{ds}=2s\dot{\gamma}$. The Euler--Lagrange equation satisfied by
$\gamma$ is%
\[
\nabla_{\dot{\gamma}}\dot{\gamma}=\frac{1}{2}\nabla R-2\operatorname*{Rc}%
(\dot{\gamma})-\frac{1}{2\tau}\dot{\gamma}.
\]
In terms of $s$, this becomes%
\[
\nabla_{\gamma^{\prime}}\gamma^{\prime}=2s^{2}\nabla R-4s\operatorname*{Rc}%
(\gamma^{\prime}),
\]
which is nonsingular at $s=0$. The computation%
\begin{align*}
\frac{d}{ds}|\gamma^{\prime}|^{2}  &  =\frac{d\tau}{ds}\frac{\partial
}{\partial\tau}g(\gamma^{\prime},\gamma^{\prime})+2g(\nabla_{\gamma^{\prime}%
}\gamma^{\prime},\gamma^{\prime})\\
&  =4s^{2}\left\langle \nabla R,\gamma^{\prime}\right\rangle
-4s\operatorname*{Rc}(\gamma^{\prime},\gamma^{\prime})
\end{align*}
shows that $|\gamma^{\prime}|$ satisfies the differential inequalities%
\begin{equation}
\frac{d}{ds}|\gamma^{\prime}|\leq2ks|\gamma^{\prime}|+2As^{2}
\label{BoundSpeedAbove}%
\end{equation}
and%
\begin{equation}
\frac{d}{ds}|\gamma^{\prime}|\geq-2Ks|\gamma^{\prime}|-2As^{2}.
\label{BoundSpeedBelow}%
\end{equation}
Let $s_{0}=\sqrt{\tau_{0}}$ and $s_{1}=\sqrt{\tau_{1}}$. Define%
\[
\overline{\psi}(s)=\left(  |\gamma^{\prime}(s_{0})|+\frac{As_{1}}{k}\right)
e^{k(s^{2}-s_{0}^{2})}-\frac{As_{1}}{k}%
\]
and%
\[
\underline{\psi}(s)=\left(  |\gamma^{\prime}(s_{0})|+\frac{As_{1}}{K}\right)
e^{K(s_{0}^{2}-s^{2})}-\frac{As_{1}}{K},
\]
replacing these by their limits if either $k$ or $K$ is zero. Note that
$\underline{\psi}(s_{0})=|\gamma^{\prime}(s_{0})|=\overline{\psi}(s_{0})$. It
is readily verified that $\overline{\psi}$ is a supersolution of
(\ref{BoundSpeedAbove}) and that $\underline{\psi}$ is a subsolution of
(\ref{BoundSpeedBelow}). So one has $\underline{\psi}(s)\leq|\frac{d\gamma
}{ds}|\leq\overline{\psi}(s)$ for $s_{0}\leq s\leq s_{1}$, as claimed.
\end{proof}

\begin{corollary}
\label{GeodesicSpeedImproved}Assume that $\operatorname*{Rc}\geq-kg$ on
$\mathcal{M}^{n}\times\lbrack0,\bar{\tau}]$. Then for any $\lambda>0$ and
$\tau^{\ast}\in(0,\bar{\tau}$), there exist positive constants $\eta$ and $C$
such that for any minimizing $\mathcal{L}$-geodesic $\gamma$ from $(\bar
{x},0)$ to $(x,\tau)\in B_{0}(\lambda)\times(0,\tau^{\ast}]$, one has%
\[
\min_{\lbrack0,\tau]}\left(  \sqrt{\sigma}|\frac{d\gamma}{d\sigma}|\right)
\geq\eta\max_{\lbrack0,\tau]}\left(  \sqrt{\sigma}|\frac{d\gamma}{d\sigma
}|\right)  -C.
\]
Furthermore, for all $\sigma\in(0,\tau]$, one has%
\[
|\frac{d\gamma}{d\sigma}|^{2}\leq\frac{2}{\eta^{2}}\left[  \frac{\ell
(\gamma(\tau),\tau)+C^{2}}{\sigma}+\frac{nk}{3}\right]  .
\]

\end{corollary}

\begin{proof}
By Lemma \ref{Neighborhood}, there exists a neighborhood $\mathcal{U}$
containing the image of $\gamma$ such that $\operatorname*{Rc}\leq Kg$ and
$|\nabla R|\leq A$ in $\mathcal{U}\times\lbrack0,\tau^{\ast}]$. Using this,
the first statement is easy to verify.

To prove the second statement, let $x=\gamma(\tau)$, so that $L(x,\tau
)=\mathcal{L}(\gamma)$. Then as in Lemma \ref{BoundsForReducedDistance}, one
has%
\[
L(x,\tau)+\frac{2nk}{3}\tau^{3/2}\geq\int_{0}^{\hat{\tau}}\sqrt{\sigma}%
|\frac{d\gamma}{d\sigma}|^{2}\,d\sigma
\]
for any $\hat{\tau}\in(0,\tau]$. Let $\psi=\min_{[0,\tau]}(\sqrt{\sigma}%
|\frac{d\gamma}{d\sigma}|)$ and $\Psi=\max_{[0,\tau]}(\sqrt{\sigma}%
|\frac{d\gamma}{d\sigma}|)$. Then for any $\delta\in(0,\hat{\tau})$ one has%
\begin{align*}
L(x,\tau)+\frac{2nk}{3}\tau^{3/2}  &  \geq2\sqrt{\hat{\tau}}\psi^{2}\\
&  \geq2\sqrt{\hat{\tau}}(\frac{\eta^{2}}{2}\Psi^{2}-C^{2})\\
&  \geq\frac{\eta^{2}}{2}\frac{\sqrt{\hat{\tau}}}{\sqrt{\hat{\tau}}-\sqrt
{\hat{\tau}-\delta}}\int_{\hat{\tau}-\delta}^{\hat{\tau}}\sqrt{\sigma}%
|\frac{d\gamma}{d\sigma}|^{2}\,d\sigma-2C^{2}\sqrt{\hat{\tau}}.
\end{align*}
Consequently, one obtains%
\begin{equation}
\int_{\hat{\tau}-\delta}^{\hat{\tau}}\sqrt{\sigma}|\frac{d\gamma}{d\sigma
}|^{2}\,d\sigma\leq\frac{\delta}{\eta^{2}\hat{\tau}}\left[  L(x,\tau
)+2C^{2}\sqrt{\hat{\tau}}+\frac{2nk}{3}\tau^{3/2}\right]  ,
\label{ShortIntegral}%
\end{equation}
whence the second statement follows.
\end{proof}

\begin{lemma}
\label{TimeDifference}Assume $\operatorname*{Rc}\geq-kg$ on $\mathcal{M}%
^{n}\times\lbrack0,\bar{\tau}]$. Let $\lambda>0$ and $\tau^{\ast}<\bar{\tau}$
be given.

\begin{enumerate}
\item There exists $C$ such that for all $x\in B_{0}(\lambda)$ and $\tau
\in(0,\tau^{\ast}]$, one has
\[
|L(x,\tau\pm\delta)-L(x,\tau)|\leq C(\frac{1}{\sqrt{\tau}}+\sqrt{\tau})\delta
\]
whenever $\delta\in(0,\tau/3)$ and $\tau\pm\delta\in\lbrack0,\tau^{\ast}]$.

\item There exists $C$ such that for all $x\in B_{0}(\lambda)$ and $\tau
\in(0,\tau^{\ast}]$, one has%
\[
|L(x,\tau\pm\delta)-L(x,\tau)|\leq C\left(  \frac{L}{\tau}+\frac{1}{\sqrt
{\tau}}+\sqrt{\tau}\right)  \delta
\]
whenever $\delta\in(0,\tau/3)$ and $\tau\pm\delta\in\lbrack0,\tau^{\ast}]$.
\end{enumerate}
\end{lemma}

\begin{proof}
Let $\alpha$ be a minimizing $\mathcal{L}$-geodesic from $(0,\bar{x})$ to
$(x,\tau)$. By Lemma \ref{Neighborhood}, we may assume that
$\operatorname*{Rc}<Kg$ and $|\nabla R|\leq A$ in $\mathcal{U}\times
\lbrack0,\tau^{\ast}]$, where $\mathcal{U}$ is a neighborhood of the image of
$\alpha$.

To bound $L$ at a later time in terms of $L$ at an earlier time, let $\beta$
denote the constant path $\beta(\sigma)=x$ for $\tau\leq\sigma\leq\tau+\delta
$. Because $\alpha$ is minimizing and $\mathcal{L}$ is additive, one has
$L(x,\tau)=\mathcal{L}(\alpha)$ and $L(x,\tau+\delta)\leq\mathcal{L}%
(\alpha)+\mathcal{L}(\beta)$. Hence there exists $C_{n}$ depending only on $n$
such that%
\[
L(x,\tau+\delta)-L(x,\tau)\leq\mathcal{L}(\beta)=\int_{\tau}^{\tau+\delta
}\sqrt{\sigma}R\,d\sigma\leq\lbrack C_{n}(k+K)\sqrt{\tau}]\delta.
\]
To bound $L$ at an earlier time in terms of $L$ at a later time, define a path
$\gamma$ from $(\bar{x},0)$ to $(x,\tau-\delta)$ by%
\[%
\begin{array}
[c]{ll}%
\gamma(\sigma)=\alpha(\sigma) & 0\leq\sigma\leq\tau-2\delta\\
\mathstrut & \mathstrut\\
\gamma(\sigma)=\alpha(2\sigma-(\tau-2\delta)) & \tau-2\delta<\sigma\leq
\tau-\delta.
\end{array}
\]
Observe that the image of $\gamma$ lies in $\mathcal{U}$ and that%
\begin{align*}
\mathcal{L}(\gamma)  &  \leq\mathcal{L}(\alpha)-\int_{\tau-2\delta}^{\tau
}\sqrt{\sigma}R(\alpha(\sigma))\,d\sigma\\
&  +4\int_{\tau-2\delta}^{\tau-\delta}\sqrt{\sigma}|\frac{d\alpha}{d\sigma
}|^{2}\,d\sigma+\int_{\tau-2\delta}^{\tau-\delta}\sqrt{\sigma}R(\gamma
(\sigma))\,d\sigma
\end{align*}
By Part (1) of Lemma \ref{GeodesicSpeed}, there exists $C^{\prime}$ such that
$|\frac{d\alpha}{d\sigma}|^{2}\leq C^{\prime}/\tau$ for $\sigma\geq
\tau-2\delta\geq\tau/3$. Since $\mathcal{L}(\alpha)=L(x,\tau)$, it follows
that%
\[
L(x,\tau-\delta)-L(x,\tau)\leq\mathcal{L}(\gamma)-\mathcal{L}(\alpha)\leq
C_{n}\left[  \frac{C^{\prime}}{\sqrt{\tau}}+(k+K)\sqrt{\tau}\right]  \delta.
\]
This proves the first statement.

To prove the second statement, use (\ref{ShortIntegral}) to estimate
$\int_{\tau-2\delta}^{\tau-\delta}\sqrt{\sigma}|\frac{d\alpha}{d\sigma}%
|^{2}\,d\sigma$.
\end{proof}

\begin{lemma}
\label{Lipschitz}If $\operatorname*{Rc}\geq-kg$ on $\mathcal{M}^{n}%
\times\lbrack0,\bar{\tau}]$, then $\ell:\mathcal{M}^{n}\times(0,\bar{\tau}%
)\ $is locally Lipschitz.

\begin{enumerate}
\item For any $\lambda>0$ and $\tau^{\ast}<\bar{\tau}$, there exists $C$ such
that%
\[
|\ell_{\tau}+\frac{\ell}{2\tau}|\leq C(\frac{1}{\tau}+1)
\]
everywhere in $B_{0}(\lambda)$ and almost everywhere in $(0,\tau^{\ast}]$, and
such that%
\[
|\nabla\ell|\leq C(\frac{1}{\tau}+1)
\]
everywhere in $(0,\tau^{\ast}]$ and almost everywhere in $B_{0}(\lambda)$.

\item There exists $C$ such that
\[
|\ell_{\tau}|\leq C\left(  \frac{\ell+1}{\tau}+1\right)
\]
everywhere in $\mathcal{M}^{n}$ and almost everywhere in $(0,\tau^{\ast}]$,
and such that%
\[
|\nabla\ell|^{2}\leq C\left(  \frac{\ell+1}{\tau}+1\right)
\]
everywhere in $(0,\tau^{\ast}]$ and almost everywhere in $\mathcal{M}^{n}$.
\end{enumerate}
\end{lemma}

\begin{proof}
We again apply Lemma \ref{Neighborhood} to get bounds $\operatorname*{Rc}\leq
Kg$ and $|\nabla R|\leq A$ on $B_{0}(2\lambda^{\ast})\times\lbrack0,\tau
^{\ast}]$, where $B_{0}(\lambda^{\ast})$ is a neighborhood of any minimizing
geodesic from $(\bar{x},0)$ to a point $(x,\tau)\in B_{0}(\lambda
)\times(0,\tau^{\ast}]$.

Wherever it is smooth, $\ell$ satisfies $\ell_{\tau}+\frac{\ell}{2\tau}%
=\frac{1}{2\sqrt{\tau}}L_{\tau}$. Thus local Lipschitz continuity in time and
the estimates for $\ell_{\tau}$ follow directly from Lemma
\ref{TimeDifference} and Rademacher's Theorem.

To show local Lipschitz continuity in space, let $x,y\in B_{0}(\lambda)$ and
$\tau\in(0,\tau^{\ast}]$ be given. We may assume that $d_{\tau}(x,y)\in
(0,\tau/3)$. Let $\alpha$ be a minimizing $\mathcal{L}$-geodesic from
$(\bar{x},0)$ to $(x,\tau)$ and let $\beta$ be a unit-speed$\ g(\tau
)$-geodesic from $x$ to $y$. Let $\delta=d_{\tau}(x,y)$ and define a path
$\gamma$ from $(\bar{x},0)$ to $(y,\tau)$ by%
\[%
\begin{array}
[c]{ll}%
\gamma(\sigma)=\alpha(\sigma) & 0\leq\sigma\leq\tau-2\delta\\
\mathstrut & \mathstrut\\
\gamma(\sigma)=\alpha(2\sigma-(\tau-2\delta)) & \tau-2\delta<\sigma\leq
\tau-\delta\\
\mathstrut & \mathstrut\\
\gamma(\sigma)=\beta(\sigma-(\tau-\delta)) & \tau-\delta<\sigma\leq\tau.
\end{array}
\]
Observe that the image of $\gamma$ belongs to $B_{0}(2\lambda^{\ast})$.
Exactly as in the proof of Lemma~\ref{TimeDifference}, one finds there exist
$C_{n}$ and $C^{\prime}$ such that%
\begin{align*}
\mathcal{L}(\gamma)  &  \leq\mathcal{L}(\alpha)-\int_{\tau-2\delta}^{\tau
}\sqrt{\sigma}R(\alpha(\sigma))\,d\sigma\\
&  +4\int_{\tau-2\delta}^{\tau-\delta}\sqrt{\sigma}|\frac{d\alpha}{d\sigma
}|^{2}\,d\sigma+\int_{\tau-\delta}^{\tau}\sqrt{\sigma}|\frac{d\beta}{d\sigma
}|^{2}\,d\sigma+\int_{\tau-2\delta}^{\tau}\sqrt{\sigma}R(\gamma(\sigma
))\,d\sigma\\
&  \leq\mathcal{L}(\alpha)+C_{n}\left[  \frac{C^{\prime}}{\sqrt{\tau}%
}+(k+K+e^{k\tau^{\ast}})\sqrt{\tau}\right]  \delta.
\end{align*}
Since $\alpha$ is minimizing, this implies that
\[
L(y,\tau)-L(x,\tau)\leq C(\frac{1}{\sqrt{\tau}}+\sqrt{\tau})d_{\tau}(x,y).
\]
Reversing the roles of $x$ and $y$ gives the same inequality for
$L(x,\tau)-L(y,\tau)$. The first gradient estimate then follows by
Rademacher's Theorem.

To prove the second gradient estimate, observe that local Lipschitz continuity
of $L$ implies that the $\mathcal{L}$-geodesic cut locus is a set of measure
zero. If $(x,\tau)$ is not in the cut locus, then the first variation formula
\cite[(7.1)]{Perelman1} implies that $\nabla L(x,\tau)=2\sqrt{\tau}%
\frac{d\alpha}{d\tau}$. The second gradient formula now follows from Corollary
\ref{GeodesicSpeedImproved}.
\end{proof}

\subsection{Integration over reduced-volume heatballs}

If $v$ is the reduced-volume density and $\varphi:\mathcal{M}^{n}\times
(0,\bar{\tau})\rightarrow\mathbb{R}$ is a given function, then the function
$P_{\varphi,v}(r)$ defined in (\ref{DefineP}) may be written as
\[
P_{\varphi,v}(r)=\int_{E_{r}}F\varphi\,d\mu\,dt,
\]
where%
\[
F=|\nabla\ell|^{2}+R(n\log\frac{r}{\sqrt{4\pi\tau}}-\ell).
\]

\begin{lemma}
\label{Finiteness}Assume that $\operatorname*{Rc}\geq-kg$ on $\mathcal{M}%
^{n}\times\lbrack0,\bar{\tau}]$. Then for any $\tau^{\ast}\in(0,\bar{\tau})$,
there exists $C$ independent of $\varphi$ such that%
\[
\frac{|P_{\varphi,v}(r)|}{r^{n}}\leq C\sup_{\mathcal{M}^{n}\times(0,cr^{2}%
)}|\varphi|
\]
whenever $0<r^{2}\leq\min\{\tau^{\ast}/c,4\pi\}$, where $c=e^{4k\bar{\tau}%
/3}/(4\pi)$.
\end{lemma}

\begin{proof}
For $0<\tau\leq\tau^{\ast}$, Part (2) of Lemma \ref{Lipschitz} implies that%
\[
|\nabla\ell|^{2}\leq\frac{C\ell+C^{\prime}}{\tau}%
\]
almost everywhere in a precompact neighborhood $\mathcal{U}$ of $\bar{x}$.
Here and in the rest of the proof, $C,C^{\prime},C^{\prime\prime}$ denote
positive constants that may change from line to line. By Corollary
\ref{BoundReducedBalls}, we may assume that $\mathcal{U}\times\lbrack
0,\tau^{\ast}]$ contains $E_{r}$ for all $r>0$ under consideration. Lemma
\ref{BoundsForReducedDistance} implies that%
\[
|\nabla\ell|^{2}\leq C\frac{d_{0}^{2}(x)}{\tau^{2}}+\frac{C^{\prime}}{\tau}%
\]
almost everywhere in $\mathcal{U}$. Let $\lambda=e^{-k\tau}\rho(r,k,\tau
)/(2\sqrt{n})$, where $\rho(r,k,\tau)$ is defined by (\ref{Define-rho}).
Together, Lemmata \ref{BoundsForReducedDistance} and \ref{rhoEstimates} show
that%
\[
0<n\log\frac{r}{\sqrt{4\pi\tau}}-\ell\leq\frac{n}{\tau}\lambda^{2}\leq\frac
{C}{\tau}(1+\tau^{2})
\]
everywhere in $E_{r}$. Hence%
\[
|F|\leq|\nabla\ell|^{2}+n(k+K)(n\log\frac{r}{\sqrt{4\pi\tau}}-\ell)\leq
C\frac{d_{0}^{2}(x)}{\tau^{2}}+\frac{C^{\prime}}{\tau}%
\]
almost everywhere in $E_{r}$. Since the volume forms $d\mu(\tau)$ are all
comparable on $B_{0}(C\lambda)\times\lbrack0,\tau^{\ast}]$, it follows from
the definition (\ref{Define-rho}) of $\rho(r,k,\tau)=2\sqrt{n}e^{k\tau}%
\lambda$ that%
\begin{align*}
\int_{B_{0}(C\lambda)}|F|\,d\mu &  \leq C^{\prime}\frac{\lambda^{n+2}}%
{\tau^{2}}+C^{\prime\prime}\frac{\lambda^{n}}{\tau}\\
&  \leq C^{\prime}\left[  \tau^{n}+\tau^{\frac{n}{2}-1}(\log\frac{r^{2}}%
{4\pi\tau})^{\frac{n}{2}+1}\right]  +C^{\prime\prime}\left[  \tau^{n-1}%
+\tau^{\frac{n}{2}-1}(\log\frac{r^{2}}{4\pi\tau})^{\frac{n}{2}}\right]  .
\end{align*}
For $r>0$ and $n\geq2$, the substitution $z=\tau/r^{2}$ shows that%
\[
\int_{0}^{cr^{2}}\tau^{\frac{n}{2}-1}(\log\frac{r^{2}}{4\pi\tau}%
)^{\frac{n+1\pm1}{2}}\,d\tau=r^{n}\int_{0}^{c}z^{\frac{n}{2}-1}(\log\frac
{1}{4\pi z})^{\frac{n+1\pm1}{2}}dz\leq Cr^{n}.
\]
Hence by Corollary \ref{BoundReducedBalls}, one has%
\[
\int_{E_{r}}|F|\,d\mu\,dt\leq\int_{0}^{cr^{2}}\left(  \int_{B_{0}(C\lambda
)}|F|\,d\mu\right)  \,d\tau\leq C^{\prime}r^{n}%
\]
whenever $0<r^{2}\leq\min\{\bar{\tau}/c,4\pi\}$. The result follows.
\end{proof}

\end{document}